\documentclass[10pt]{article}
\usepackage{amsmath, amsfonts, amsthm,amssymb}
\usepackage{graphicx}
\begin{document}
\title{
Change of variable formulas for \\
non-anticipative functionals on path space}
\author{Rama CONT \footnote{Laboratoire de Probabilit\'es et Mod\`eles Al\'eatoires, UMR 7599 CNRS-Universit\'e de Paris VI, France \& Columbia University, New York. Email: Rama.Cont@columbia.edu}
\hskip 2cm \hfill  David-Antoine FOURNIE \footnote{Dept of
Mathematics, Columbia University, New York. Email:
df2243@columbia.edu}
}
\date{\today\footnote{We thank  Jean Jacod and Paul Malliavin for helpful discussions.}}
\newtheorem{theorem}{Theorem}
\newtheorem{corollary}[theorem]{Corollary}
\newtheorem{lemma}[theorem]{Lemma}
\newtheorem{proposition}[theorem]{Proposition}

\theoremstyle{remark}
\newtheorem{remark}{Remark}
\newtheorem{example}{Example}

\theoremstyle{definition}
\newtheorem{definition}{Definition}
\newtheorem{assumption}{Assumption}
\newcommand{\Ut}{\mathcal{U}_t}
\newcommand{\Et}{\mathcal{U}_t\times \mathcal{S}_t}
\newcommand{\pp}{\mathbb{P}}
\newcommand{\qq}{\mathbb{Q}}
\newcommand{\ba}{\begin{eqnarray}}
\newcommand{\ea}{\end{eqnarray}}

\maketitle

\begin{abstract}
We  derive a functional change of variable formula for {\it
non-anticipative} functionals defined on the space of
$\mathbb{R}^d$-valued right continuous paths with left limits. The
functional is only required to possess  certain directional
derivatives, which may be computed pathwise. Our results lead to
functional extensions of the Ito formula for a large class of
stochastic processes, including semimartingales and Dirichlet
processes. In particular, we show the stability of the  class of
semimartingales under certain functional transformations.
\end{abstract}
\newpage\tableofcontents \newpage
In his seminal paper {\it Calcul d'Ito sans probabilit\'es}
\cite{follmer79}, Hans F\"ollmer proposed a non-probabilistic
version of the Ito formula \cite{ito44}: F\"ollmer showed that if a
real-valued cadlag (right continuous with left limits) function $x$
has  finite quadratic variation  along a sequence
$\pi_n=(t^{n}_k)_{k=0..n}$ of subdivisions  of  $[0,T]$ with  step
size decreasing to zero, in the sense that the sequence of discrete
measures
$$ \sum_{k=0}^{n-1} \|x(t_{k+1}^n)-x(t^n_k)\|^2\delta_{t^n_k}$$
converges vaguely to a Radon measure with Lebesgue decomposition
$\xi+ \sum_{t\in[0,T]} |\Delta x(t)|^2 \delta_t$ then for  $f\in
C^1(\mathbb{R})$ one can define the pathwise integral \ba \int_0^T
f(x(t)) d^\pi x =  \lim_{n\to\infty} \sum_{i=0}^{n-1} f(x({t^n_i}))
.(x(t^n_{i+1})-x(t^n_i)) \label{follmerintegral.def}\ea as a limit
of Riemann sums along the subdivision $\pi=(\pi_n)_{n\geq 1}$. In
particular if $X=(X_t)_{t\in[0,T]}$ is a semimartingale
\cite{dm,meyer76,protter}, which is the classical setting for
stochastic calculus, the paths of $X$ have finite quadratic
variation along such subsequences: when applied to the paths of $X$,
F\"ollmer's integral \eqref{follmerintegral.def} then coincides,
with probability one, with the Ito stochastic integral $\int_0^T
f(X) dX$ with respect to the semimartingale $X$.
  This construction may in fact be carried out for a more general class of processes, including the class of Dirichlet processes \cite{coquet03,follmer79,follmer81,lyons94}.

  Of course, the Ito stochastic integral with respect to a semimartingale $X$ may be defined  for a much larger class of integrands: in particular, for a caglad process $Y$ defined as a  {\it non-anticipative functional } $Y(t)=F_t(X(u),0\leq u\leq t)$ of $X$, the stochastic integral $\int_0^T Y dX$ may be defined as a limit of non-anticipative Riemann sums \cite{protter}.

   Using a notion of directional derivative for functionals proposed by Dupire \cite{dupire09}, we
   extend F\"ollmer's pathwise
      change of variable formula to non-anticipative functionals on the space $D([0,T],\mathbb{R}^d)$
      of cadlag paths (Theorem \ref{itonoprobability.theorem}).
The requirement on the functionals is to possess certain directional
derivatives which may be computed pathwise.
    Our construction allows to define a pathwise integral $\int F_t(x)dx$, defined as a limit of Riemann sums, for a class of functionals $F$
    of a cadlag path
    $x$ with finite quadratic variation.
     Our results
lead to functional extensions of the Ito formula for semimartingales
(Section \ref{semimartingale.sec}) and Dirichlet processes (Section
\ref{dirichlet.sec}).  In particular, we show the stability of the
the class of semimartingales under functional transformations
verifying a regularity condition. These results yield a
non-probabilistic proof for functional Ito formulas obtained in
\cite{ContFournie09a,ContFournie09b,dupire09} using probabilistic
methods and extend them to the case of discontinuous
semimartingales.

{\bf Notation}\\ \noindent   For a path $x\in
D([0,T],\mathbb{R}^d)$, denote by $x(t)$ the value of $x$ at $t$ and
by $x_t=(x(u), 0\leq u\leq t)$ the restriction of $x$ to $[0,t]$.
Thus  $x_t\in D([0,t],\mathbb{R}^d)$. For a stochastic process $X$
we shall similarly denote $X(t)$ its value at $t$ and
 $X_t=(X(u), 0\leq u\leq t)$ its path on $[0,t]$.

\section{Non-anticipative functionals on spaces of paths}

Let $T>0$, and $U\subset \mathbb{R}^d$ be an open subset of
$\mathbb{R}^d$ and $S \subset \mathbb{R}^m$ be a Borel subset of
$\mathbb{R}^m$.
 We call "$U$-valued cadlag function" a right-continuous function $f:[0,T]\mapsto U$ with left limits such that
 for each $t\in ]0,T]$, $f(t-)\in U$. Denote by $\mathcal{U}_t=D([0,t],U)$ (resp.
$\mathcal{S}_t=D([0,t],S)$ the space of $U$-valued cadlag functions
(resp. $S$), and $C_0([0,t],U)$ the set of continuous functions with
values in $U$.

When dealing with functionals of a path $x(t)$  indexed by time, an
important class is formed by those which are {\it non-anticipative},
in the sense that they only depend on the past values of $x$.
 A family $Y: [0,T]\times \mathcal{U}_T\mapsto \mathbb{R}$ of functionals
is said to be {\it non-anticipative} if, for all $(t,x)\in
[0,T]\times \mathcal{U}_T$, $Y(t,x)=Y(t,x_t)$ where $x_t=x_{|[0,t]}$
denotes the restriction of the path $x$ to $[0,t]$. A
non-anticipative functional may thus be represented as
$Y(t,x)=F_t(x_t)$ where $(F_t)_{t\in[0,T]}$ is a family of maps
$F_t:\mathcal{U}_t\mapsto \mathbb{R}$. 
This motivates the following definition:
\begin{definition}[Non-anticipative  functionals on path space]
A non-anticipative functional on $\mathcal{U}_T$  is a family
$F=(F_t)_{t \in [0,T]}$  of maps \ba F_t: \mathcal{U}_t \to
\mathbb{R}\nonumber\ea
\end{definition}
 $Y$ is said to be {\it predictable}\footnote{This notion coincides
with the usual definition of predictable process when the path space
$\mathcal{U}_T$ is endowed with the filtration of the canonical
process, see Dellacherie \& Meyer \cite[Vol. I]{dm}.} if, for all
$(t,x)\in [0,T]\times \mathcal{U}_T$, $Y(t,x)=Y(t,x_{t-})$ where
$x_{t-}$ denotes the function defined on $[0,t]$ by
$$ x_{t-}(u)=x(u)\quad u\in[0,t[\qquad x_{t-}(t)=x(t-)$$
Typical examples of predictable functionals are  integral
functionals, e.g.
$$ Y(t,x)= \int_0^t G_s(x_s) ds$$
where $G$ is a non-anticipative, locally integrable, functional.

 If $Y$ is predictable then $Y$ is non-anticipative, but
predictability is a stronger property. Note that $x_{t-}$ is cadlag
and should \textit{not} be confused with the caglad path $u\mapsto
x(u-)$.

We consider throughout this paper non-anticipative functionals  \ba
F=(F_t)_{t \in [0,T]}\qquad F_t: \mathcal{U}_t \times
\mathcal{S}_t\to \mathbb{R}\nonumber\ea where $F$ has a predictable
dependence with respect to the second argument: \ba \forall t \leq
T,\quad \forall (x,v) \in \mathcal{U}_t \times \mathcal{S}_t,\quad
F_t(x_t,v_t)=F_t(x_t,v_{t-})\label{predictable.eq} \ea $F$ can be
viewed as a functional on the vector bundle $\Upsilon=\bigcup_{t \in
[0,T]} \mathcal{U}_t \times \mathcal {S}_t$. We will also consider
non-anticipative functionals $F=(F_t)_{t \in [0,T[}$ indexed by
$[0,T[$.

\subsection{Horizontal and vertical perturbation of a path}\label{extensions.sec}

Consider a path $x\in D([0,T]),U)$ and denote by $x_t\in
\mathcal{U}_t$ its restriction to $[0,t]$ for $t<T$. For
$h\geq 0$, the {\it horizontal} extension $x_{t,h}\in
D([0,t+h],\mathbb{R}^d)$ of $x_t$ to $[0,t+h]$ is defined as\ba
x_{t,h}(u)=x(u) \qquad u\in[0,t]\ ; \qquad x_{t,h}(u)=x(t) \qquad
u\in ]t,t+h]\ea For $h\in\mathbb{R}^d$ small enough, we define the {\it vertical}
perturbation $x^h_t$ of $x_t$  as the cadlag path obtained by
shifting the endpoint by $h$: \ba x^h_t(u)=x_t(u) \qquad u\in[0,t[ &
x^h_t(t)=x(t)+h\ea or in other words $x^h_t(u)=x_t(u)+h 1_{t=u}$. By convention, $x^u_{t,h}=(x^u_t)_{t,h}$, ie the vertical perturbation precedes the horizontal extension.\\

We now define a distance between two paths, not necessarily defined
on the same time interval. For $T\geq t'=t+h\geq t\geq 0$, $(x,v)\in
\mathcal{U}_t\times {S}^+_t$ and $(x',v') \in
D([0,t+h],\mathbb{R}^d)\times \mathcal{S}_{t+h}$ define \ba
d_{\infty}(\ (x,v),(x',v')\ )=
\sup_{u\in[0,t+h]}|x_{t,h}(u)-x'(u)|+\sup_{u\in[0,t+h]}|v_{t,h}(u)-v'(u)|+\
h \ea If the paths $(x,v),(x',v')$ are defined on the same time
interval, then $d_\infty((x,v),(x',v'))$ is simply the distance in
supremum norm.

\subsection{Classes of non-anticipative functionals}\label{regularity.sec}

Using the  distance $d_\infty$ defined above, we now introduce
various notions of continuity for non-anticipative functionals.
\begin{definition}[Continuity at   fixed times]
A non-anticipative functional $F=(F_t)_{t\in[0,T]}$ is said to be
continuous at fixed times if and only if for any $t \leq T$, $F_t:
\mathcal{U}_t \times \mathcal{S}_t\mapsto \mathbb{R}$ is continuous
for the supremum norm.
\end{definition}

\begin{definition}[Left-continuous functionals]
Define  $\mathbb{F}^\infty_l$ as the set of functionals
$F=(F_t,t\in[0,T])$  which verify: \ba \forall t \in
[0,T],\quad\forall \epsilon > 0, \forall (x,v)\in \mathcal{U}_t
\times \mathcal{S}_t, \quad \exists \eta > 0, \forall h \in [0,t],
\quad \nonumber\\ \forall (x',v') \in \mathcal{U}_{t-h} \times
\mathcal{S}_{t-h},\qquad d_\infty((x,v),(x',v')) < \eta \Rightarrow
|F_t(x,v)-F_{t-h}(x',v')| < \epsilon \ea \label{lcontinuous.def}
\end{definition}

\begin{definition}[Right-continuous functionals]
Define  $\mathbb{F}^\infty_r$ as the set of functionals
$F=(F_t,t\in[0,T[)$  which verify \ba \forall t \in
[0,T],\quad\forall \epsilon > 0, \forall (x,v)\in \mathcal{U}_t
\times \mathcal{S}_t, \quad \exists \eta > 0, \forall  h \in
[0,T-t],  \nonumber\\ \forall (x',v') \in \mathcal{U}_{t+h} \times
\mathcal{S}_{t+h},\qquad d_\infty((x,v),(x',v')) < \eta \Rightarrow
|F_t(x,v)-F_{t+h}(x',v')| < \epsilon \ea \label{rcontinuous.def}
\end{definition}
We denote $\mathbb{F}^\infty= \mathbb{F}^\infty_r\cap
\mathbb{F}^\infty_r$ the set of continuous non-anticipative
functionals.

We  call a functional "boundedness preserving" if it is
bounded on  each bounded set of  paths:
\begin{definition}[ Boundedness-preserving functionals]
Define $\mathbb{B}$ as the set of non-anticipative functionals $F$
such that for every compact subset $K$ of $U$, every $R >0$, there
exists a constant $C_{K,R}$ such that: \ba \forall t \leq T, \forall
(x,v) \in D([0,t],K) \times \mathcal{S}_t,\sup_{s \in [0,t]} |v(s)|
< R \Rightarrow |F_t(x,v)| < C_{K,R}\label{boundedpreserving.eq} \ea
\label{boundedness.def}
\end{definition}

In particular if $F\in \mathbb{B}$, it is "locally" bounded in the neighborhood of any given path i.e.
\ba \forall (x,v) \in \mathcal{U}_T \times \mathcal{S}_T,\quad \exists C>0, \eta >0, \qquad
\forall t \in [0,T],\quad \forall (x',v') \in \mathcal{U}_t \times \mathcal{S}_t,\nonumber\\
d_\infty((x_t,v_t),(x',v')) < \eta \Rightarrow \forall t \in [0,T], |F_t(x',v')| \leq C  \label{localbounded.eq} \ea

The following result describes the behavior of paths generated by the functionals in the above classes:
\begin{proposition}[Pathwise regularity]\ \\
\begin{enumerate}
\item
If $F\in \mathbb{F}^{\infty}_l$ then for any $(x,v)\in
\mathcal{U}_T \times \mathcal{S}_T$, the path $t\mapsto
F_t(x_{t-},v_{t-})$ is left-continuous.
\item
If $F\in \mathbb{F}^{\infty}_r$ then for any $(x,v)\in
\mathcal{U}_T \times \mathcal{S}_T$, the path $t\mapsto
F_t(x_{t},v_{t})$ is right-continuous.
\item
If $F\in \mathbb{F}^{\infty}$ then for any $(x,v)\in
\mathcal{U}_T \times \mathcal{S}_T$, the path $t\mapsto
F_t(x_t,v_t)$ is cadlag and continuous at all points where $x$ and $v$ are continuous.
\item
If $F\in \mathbb{F}^{\infty}$ further verifies
(\ref{predictable.eq}) then for any $(x,v)\in \mathcal{U}_T \times
\mathcal{S}_T$, the path $t\mapsto F_t(x_t,v_t)$ is cadlag and
continuous at all points where $x$ is continuous.
\label{cadlag.prop}
\item If $F\in \mathbb{B}$, then for any $(x,v)\in
\mathcal{U}_T \times \mathcal{S}_T$, the path $t\mapsto
F_t(x_t,v_t)$ is bounded.
\end{enumerate}
\end{proposition}
\begin{proof}
\begin{enumerate}
\item
Let $F\in \mathbb{F}^{\infty}_l$ and $t \in [0,T)$. For $h> 0$ sufficiently small,
\begin{equation}
d_{\infty}((x_{t-h},v_{t-h}),(x_{t-},v_{t-}))=\sup_{u \in (t-h,t)} |x(u)-x(t-)| + \sup_{u \in (t-h,t)} |v(u)-v(t-)| + h
\end{equation}
Since $x$ and $v$ are cadlag, this quantity converges to 0 as $h \rightarrow 0+$, so
$$ F_{t-h}(x_{t-h},v_{t-h})-F_t(x_{t-},v_{t-}) \mathop{\to}^{h\to 0^+} 0 $$
so $t\mapsto
F_t(x_{t-},v_{t-})$ is left-continuous.
\item
Let $F\in \mathbb{F}^{\infty}_r$ and $t \in [0,T)$. For $h> 0$ sufficiently small,
\begin{equation}
d_{\infty}((x_{t+h},v_{t+h}),(x_{t},v_{t}))=\sup_{u \in [t,t+h)} |x(u)-x(t)| + \sup_{u \in [t,t+h)} |v(u)-v(t)| + h
\end{equation}
Since $x$ and $v$ are cadlag, this quantity converges to 0 as $h \rightarrow 0+$, so
$$ F_{t+h}(x_{t+h},v_{t+h})-F_t(x_{t},v_{t}) \mathop{\to}^{h\to 0^+} 0 $$
so $t\mapsto
F_t(x_{t},v_{t})$ is right-continuous.
\item
Assume now that $F$ is in $\mathbb{F}^{\infty}$ and let $t \in
]0,T]$. Denote $(\Delta x(t),\Delta v(t))$ the jump of $(x,v)$ at
time $t$. Then
\begin{equation*}
d_\infty((x_{t-h},v_{t-h}),x_t^{-\Delta x(t)},v_t^{-\Delta
v(t)}))=\sup_{u \in [t-h,t)} |x(u)-x(t)| + \sup_{u \in [t-h,t)}
|v(u)-v(t)| + h
\end{equation*}
and this quantity goes to $0$ because $x$ and $v$ have left limits.
Hence the path has left limit $F_{t}(x^{-\Delta x(t)}_t,v_t^{-\Delta
v(t)})$ at $t$. A similar reasoning proves that it has right-limit
$F_t(x_t,v_t)$.
\item
If $F\in \mathbb{F}^{\infty}$ verifies \eqref{predictable.eq}, for
 $t \in ]0,T]$ the  path $t\mapsto F_t(x_t,v_t)$ has left-limit
$F_{t}(x^{-\Delta x(t)}_t,v_t^{-\Delta v(t)})$ at $t$, but
\eqref{predictable.eq} implied that this left-limit equals
$F_t(x^{-\Delta x(t)}_t,v_t)$.
\end{enumerate}
\end{proof}

\subsection{Measurability properties}
Consider, on the path space $\mathcal{U}_T \times \mathcal{S}_T$,
the filtration $(\mathcal{F}_t)$ generated by the canonical process
\ba (X,V): \mathcal{U}_T\times\mathcal{S}_T \times [0,T]&\mapsto &
U\times S\nonumber\\
(x,v),t&\to & (X,V)((x,v),t)=(x(t),v(t)) \ea ${\mathcal F}_t$ is the
smallest sigma-algebra on $\mathcal{U}_T\times\mathcal{S}_T$ such
that all coordinate maps $(X(.,s),V(.,s)), s\in[0,t]$ are ${\mathcal
F}_t$-measurable.

 The following result,
proved in Appendix \ref{measurabilityproof.sec}, clarifies the
measurability properties of processes defined by functionals in
$\mathbb{F}^{\infty}_l,\mathbb{F}^{\infty}_r$:
\begin{theorem}\label{Measurability.thm}
If $F$ is continuous at fixed time, then the process $Y$ defined by
$Y((x,v),t)=F_t(x_{t},v_{t})$ is ${\cal F}_t$-adapted. If  $F \in
\mathbb{F}^{\infty}_l$ or $F \in \mathbb{F}^{\infty}_r$, then:
\begin{enumerate}
\item the process $Y$ defined by $Y((x,v),t)=F_t(x_{t},v_{t})$ is optional.
\item the process $Z$ defined by $Z((x,v),t)=F_t(x_{t-},v_{t-})$ is predictable.
\end{enumerate}
\end{theorem}

\section{Pathwise derivatives of non-anticipative functionals}
\label{derivative.sec}

\subsection{Horizontal   derivative}
We now define a pathwise derivative for a non-anticipative functional
$F=(F_t)_{t\in[0,T]}$, which  may be seen as a ``Lagrangian"
derivative along the path $x$.
\begin{definition}[Horizontal derivative]
The {\it horizontal derivative}  at $(x,v)\in
\mathcal{U}_t\times \mathcal{S}_t$ of non-anticipative
 functional $F=(F_t)_{t\in[0,T[}$ is defined as
\ba \mathcal{D}_tF(x,v) = \lim_{h\to 0^+} \frac{F_{t+h}(
x_{t,h},v_{t,h})-F_t(x,v)}{h}\label{horizontalderivative.eq}\ea if
the corresponding limit exists. If \eqref{horizontalderivative.eq} is defined for all
$(x,v)\in\Upsilon$ the map \ba \mathcal{D}_tF :
\mathcal{U}_t\times
\mathcal{S}_t &\mapsto & \mathbb{R}^d \nonumber\\
(x,v)&\to& \mathcal{D}_tF(x,v)\ea defines a non-anticipative
functional $\mathcal{D}F=(\mathcal{D}_tF)_{t\in[0,T[}$,  the {\it
horizontal derivative} of $F$. \label{horizontalderivative.def}
\end{definition}
We will occasionally use the following ``local Lipschitz property"
that is weaker than horizontal differentiability:
\begin{definition} \label{timelipschitz.def}
A non-anticipative functional $F$ is said to have the horizontal
local Lipschitz property if and only if:
\ba \forall (x,v) \in \mathcal{U}_T \times \mathcal{S}_T, \exists C>0,\eta>0, \forall t_1 < t_2 \leq T, \forall(x',v') \in \mathcal{U}_{t_1} \times \mathcal{S}_{t_1}, \nonumber \\
d_{\infty}((x_{t_1},v_{t_1}),(x',v'))< \eta \Rightarrow
|F_{t_2}(x'_{t_1,t_2-t_1},v'_{t_1,t_2-t_1})-F_{t_1}((x'_{t_1},v'_{t_1}))|<C(t_2-t_1)\label{localLipschitz.eq}
\ea
\end{definition}
\subsection{Vertical derivative}
Dupire \cite{dupire09} introduced a pathwise spatial derivative for
non-anticipative functionals, which we now introduce. Denote
$(e_i,i=1..d)$ the canonical basis in $\mathbb{R}^d$.
\begin{definition} A non-anticipative
 functional $F=(F_t)_{t\in[0,T]}$  is said to be
{\it vertically differentiable}   at $(x,v)\in
D([0,t]),\mathbb{R}^d)\times D([0,t],S^+_d)$ if \ba \mathbb{R}^d&\mapsto&\mathbb{R}\nonumber\\
e&\to & F_t(x^{e}_t,v_t) \nonumber\ea is differentiable at $0$.
 Its gradient at $0$ \ba {\nabla}_xF_t\ (x, v) = (\partial_iF_t(x,v),\
i=1..d) \qquad {\rm where}\quad
\partial_iF_t(x,v)=\lim_{h\to 0}\frac{F_t(x^{he_i}_t,v)-F_t(x,v)}{h}\quad\label{verticalderivative.eq}\ea
is called  the {\it vertical derivative} of $F_t$ at $(x,v)$. If
\eqref{verticalderivative.eq} is defined for all $(x,v)\in\Upsilon$,
the {\it vertical derivative}
   \ba {\nabla}_xF : \mathcal{U}_t\times
\mathcal{S}_t &\mapsto & \mathbb{R}^d \nonumber\\
(x,v)&\to& {\nabla}_xF_t(x,v)\ea define  a non-anticipative
functional ${\nabla}_xF=({\nabla}_xF_t)_{t\in[0,T]}$ with values in
$\mathbb{R}^d$. \label{verticalderivative.def}\end{definition}
\begin{remark}
If a vertically differentiable functional verifies
\eqref{predictable.eq}, its vertical derivative also verifies
\eqref{predictable.eq}.
\end{remark}
\begin{remark}$\partial_iF_t(x,v)$ is simply the directional derivative
 of $F_t$
in direction  $(1_{\{t\}} e_i,0)$. Note that this involves examining cadlag perturbations of the path $x$, even if $x$ is continuous.
\end{remark}
\begin{remark}
If $F_t(x,v)=f(t,x(t))$ with $f\in C^{1,1}([0,T[\times
\mathbb{R}^d)$ then we retrieve the usual partial derivatives:
$$ {\mathcal D}_tF(x,v)= \partial_tf(t,x(t))\qquad \nabla_xF_t(x_t,v_t)= \nabla_xf(t,x(t)).$$\end{remark}
\begin{remark} Note that the assumption \eqref{predictable.eq} that
$F$ is predictable with respect to the second variables entails that
for any $t\in[0,T]$, $F_t(x_t,v_t^e)=F_t(x_t,v_t)$ so an analogous
notion of derivative with respect to $v$ would be identically zero
under assumption \eqref{predictable.eq}.
\end{remark}
If $F$ admits a horizontal (resp. vertical) derivative ${\cal D}F$
(resp. $\nabla_xF$) we may iterate the operations described above
and define higher order horizontal and vertical derivatives.
\begin{definition} \label{Cab.def}
Define $\mathbb{C}^{j,k}$ as the set of functionals $F$ which are
\begin{itemize} \item continuous at fixed times, \item admit $j$ horizontal derivatives and $k$ vertical derivatives
 at all $(x,v) \in \mathcal{U}_t \times \mathcal{S}_t$, $t \in [0,T[$
  \item $\mathcal{D}^mF,m \leq
j,\nabla^n_xF, n \leq k$ are  continuous at fixed
times.\end{itemize}
\end{definition}

\section{Change of variable formula for functionals of a continuous path}
We now state our first main result, a functional change of variable
formula which extends the It\^o formula without probability due to
F{\"o}llmer \cite{follmer79} to functionals. We denote here $S^+_d$
the set of positive symmetric $d \times d$ matrices.

\begin{definition} \label{finitequadraticvariation.def}
Let $\Pi_n=(t^n_0,\ldots,t^n_{k(n)})$, where $0=t^n_0 \leq t^n_1
\leq \ldots \leq t^n_{k(n)}=T$, be a sequence of subdivisions of
$[0,T]$ with step decreasing to 0 as $n \rightarrow \infty$. $f \in
C_0([0,T],\mathbb{R})$ is said to have finite quadratic variation
along $(\pi_n)$ if the sequence of discrete measures: \ba
\xi^n=\sum_{i=0}^{k(n)-1} (f(t^n_{i+1})-f(t^n_i))^2 \delta_{t^n_i}
\ea where $\delta_t$ is the Dirac measure at $t$, converge vaguely
to a Radon measure $\xi$ on $[0,T]$ whose atomic part is null. The
increasing function $[f]$ defined by
$$[f](t)= \xi([0,t])$$
is then called the quadratic variation of   $f$ along the sequence $(\pi_n)$.\\
$x \in C_0([0,T],U)$ is said to have finite quadratic variation
along the sequence $(\pi_n)$ if the functions $x_i,1 \leq i \leq d$
and $x_i+x_j, 1 \leq i < j \leq d$ do. The quadratic variation of
$x$ along $(\pi_n)$ is the $S^+_d$-valued function $x$ defined by:
\ba [x]_{ii}=[x_i], [x]_{ij}=\frac{1}{2} ([x_i+x_j]-[x_i]-[x_j]), i
\neq j \ea

\end{definition}

\begin{theorem}[Change of variable formula for functionals of continuous paths]\label{itonoprobability.theorem}
Let $(x,v) \in C_0([0,T],U) \times \mathcal{S}_T$   such that $x$
has finite quadratic variation along $(\pi_n)$ and verifies $
\sup_{t \in [0,T]-\pi_n} |v(t)-v(t-)| \rightarrow 0$. Denote:
\ba x^n(t)=\sum_{i=0}^{k(n)-1} x(t_{i+1})1_{[t_i,t_{i+1}[}(t)++x(T)1_{\{T\}}(t)  \nonumber \\
v^n(t)=\sum_{i=0}^{k(n)-1}
v(t_{i})1_{[t_i,t_{i+1}[}(t)+v(T)1_{\{T\}}(t), \qquad
h^n_i=t^n_{i+1}-t^n_i \ea Then for any     non-anticipative
functional $F \in \mathbb{C}^{1,2}$ satisfying the following
assumptions:
\begin{enumerate}
\item $F,\nabla_xF,\nabla^2_xF \in \mathbb{F}^{\infty}_l$
\item $\nabla^2_x F,\mathcal{D}F$ satisfy the local boundedness  property  \eqref{localbounded.eq}
\end{enumerate}  the following limit
 \ba \lim_{n\to\infty}  \sum_{i=0}^{k(n)-1} \nabla_xF_{t^n_{i}}(x^n_{t^n_i-},v^n_{t^n_i-}) (x(t^n_{i+1})-x(t^n_i))\label{follmerintegral.eq}\ea
exists. Denoting this limit by $\int_0^T \nabla_xF(x_u,v_u)
d^{\pi}x$ we have
\begin{eqnarray}
F_T(x_T,v_T)-F_0(x_0,v_0) = \int_{0}^{T}
\mathcal{D}_tF_t(x_{u},v_{u}) du+ \int_{0}^{T} \frac{1}{2} {\rm
tr}\left({}^t\nabla^2_xF_t(x_{u},v_{u}) d[x](u)\right) + \int_0^T
\nabla_xF(x_u,v_u) d^{\pi}x \label{functional.noproba.ito.eq}
\end{eqnarray}
\end{theorem}
\begin{remark}[F\"ollmer integral] The limit \eqref{follmerintegral.eq}, which we call the {\it F\"ollmer integral},
was defined in \cite{follmer79} for integrands of the form $f(X(t))$
where $f\in C^1(\mathbb{R}^d)$. It depends a priori on the sequence
$\pi$ of subdivisions, hence the notation $\int_0^T
\nabla_xF(x_u,v_u) d^{\pi}x$. We will see in Section
\ref{semimartingale.sec} that when $x$ is the sample path of a {\it
semimartingale}, the limit is in fact almost-surely independent of
the choice of $\pi$.
\end{remark}
\begin{remark}
The regularity conditions on $F$ are given independently of $(x,v)$
and of the sequence of subdivisions $(\pi_n)$.
\end{remark}

\begin{proof}
Denote $\delta x^n_i=x(t^n_{i+1})-x(t^n_i)$. Since $x$ is continuous
hence uniformly continuous on $[0,T]$, and  using Lemma
\ref{uniform.lemma} for $v$, the quantity \ba
\eta_n=\sup\{|v(u)-v(t^n_i)|+|x(u)-x(t^n_i)|+|t^n_{i+1}-t^n_i|,0
\leq i \leq k(n)-1,u\in[t^n_{i},t^n_{i+1})\} \ea  converges to 0 as
$n \rightarrow \infty$. Since $\nabla^2_x F,\mathcal{D}F$ satisfy
the local boundedness  property  \eqref{localbounded.eq}, for $n$
sufficiently large there exists $C > 0$ such that $$\forall t < T,
\forall (x',v') \in \mathcal{U}_t \times \mathcal{S}_t,\qquad
d_\infty((x_t,v_t),(x',v')) < \eta_n \Rightarrow
|\mathcal{D}_tF_t(x',v')| \leq C,|\nabla^2_xF_t(x',v')| \leq C$$
Denoting $K=\overline{\{x(u), s \leq u \leq t \}}$ which is a
compact subset of $U$, and $U^{c}=\mathbb{R^d}-U$ its complement,
one can also assume $n$ sufficiently large so that $d(K,U^{c}) >
\eta_n$.

For $i \leq k(n)-1$, consider the decomposition:
\ba F_{t^n_{i+1}}(x^n_{t^n_{i+1}-},v^n_{t^n_{i+1}-})-F_{t^n_i}(x^n_{t^n_i-},v^n_{t^n_i-})&=&F_{t^n_{i+1}}(x^n_{t^n_{i+1}-},v^n_{t^n_i,h^n_i})-F_{t^n_i}(x^n_{t^n_i},v^n_{t^n_i}) \nonumber \\
&+&F_{t^n_i}(x^n_{t^n_{i}},v^n_{t^n_i-})-F_{t^n_{i}}(x^n_{t^n_i-},v^n_{t^n_i-})
 \label{24.eq}\ea

where we have used   property \eqref{predictable.eq} to have
$F_{t^n_i}(x^n_{t^n_i},v^n_{t^n_i})=F_{t^n_i}(x^n_{t^n_{i}},v^n_{t^n_i-})$.
The first term can be written $ \psi(h^n_i)-\psi(0) $ where: \ba
\psi(u)=F_{t^n_i+u}(x^n_{t^n_i,u},v^n_{t^n_i,u}) \ea Since $F \in
\mathbb{C}^{1,2}([0,T])$, $\psi$ is right-differentiable, and
moreover by lemma \ref{cadlag.prop}, $\psi$ is left-continuous, so:
\ba
F_{t^n_{i+1}}(x^n_{t^n_i,h^n_i},v^n_{t^n_i,h^n_i})-F_{t^n_i}(x^n_{t^n_i},v^n_{t^n_i})=
\int_{0}^{t^n_{i+1}-t^n_i}
\mathcal{D}_{t^n_i+u}F(x^n_{t^n_i,u},v^n_{t^n_i,u}) du \ea

The second term can be written $\phi(\delta x^n_i)-\phi(0)$, where:
\ba \phi(u)=F_{t^n_{i}}(x^{n,u}_{t^n_i-},v^n_{t^n_i-}) \ea
Since $F \in \mathbb{C}^{1,2}([0,T])$, $\phi$ is well-defined and $C^2$ on the convex set $B(x(t^n_i),\eta_n) \subset U$, with:
\ba \phi'(u)=\nabla_x F_{t^n_{i}}(x^{n,u}_{t^n_i-},v^n_{t^n_i-})  \nonumber \\
\phi''(u)=\nabla^2_x F_{t^n_{i}}(x^{n,u}_{t^n_i-},v^n_{t^n_i-}) \ea
So a second order Taylor expansion of $\phi$ at $u=0$ yields:
\ba F_{t^n_{i}}(x^n_{t^n_{i}},v^n_{t^n_i-})-F_{t^n_{i}}(x^n_{t^n_i-},v^n_{t^n_i-})=\nabla_x F_{t^n_{i}}(x^n_{t^n_i-},v^n_{t^n_i-})\delta x^n_i \nonumber \\
+ \frac{1}{2} {\rm tr}\left(\nabla^2_x
F_{t^n_{i}}(x^n_{t^n_i-},v^n_{t^n_i-})\quad {}^t\delta x^n_i \delta
x^n_i\right) +r^n_i \ea where $r^n_i$ is bounded by \ba K |\delta
x^n_i|^2 \sup_{x \in B(x(t^n_i),\eta_n)} |\nabla^2_x
F_{t^n_{i}}(x^{n,x-x(t^n_i)}_{t^n_i-},v^n_{t^n_i-})-\nabla^2_x
F_{t^n_{i}}(x^n_{t^n_i-},v^n_{t^n_i-})|  \ea Denote $i^n(t)$ the
index such that $t \in [t^n_{i^n(t)},t^n_{i^n(t)+1})$. We now sum
all the terms above from $i=0$ to $k(n)-1$:.
\begin{itemize}
\item The left-hand side of \eqref{24.eq} yields $F_{T}(x^n_{T-},v^n_{T-})-F_{0}(x_0,v_0)$, which converges to $F_{T}(x_{T-},v_{T-})-F_0(x_0,v_0)$ by left-continuity of $F$, and this quantity equals $F_{T}(x_{T},v_{T})-F_0(x_0,v_0)$ since $x$ is continuous and $F$ is predictable in the second variable.

\item The first line in the right-hand side can be written:
\ba \int_{0}^{T}
\mathcal{D}_uF(x^n_{t^n_{i^n(u)},u-t^n_{i^n(u)}},v^n_{t^n_{i^n(u)},u-t^n_{i^n(u)}})du\label{32.eq}
\ea where the integrand converges to $\mathcal{D}_uF(x_u,v_{u-})$
and is bounded by $C$. Hence the dominated  convergence theorem
applies and \eqref{32.eq} converges to: \ba \int_{0}^{T}
\mathcal{D}_uF(x_{u},v_{u-}) du = \int_{0}^{T}
\mathcal{D}_uF(x_{u},v_{u}) \ea since $v_u=v_{u-},du$-almost
everywhere.

\item The second line can be written:
\ba \sum_{i=0}^{k(n)-1} \nabla_xF_{t^n_{i}}(x^n_{t^n_i-},v^n_{t^n_i-}) (x_{t^n_{i+1}}-x_{t^n_i})
    + \sum_{i=0}^{k(n)-1} \frac{1}{2} tr[\nabla^2_x F_{t^n_{i}}(x^n_{t^n_i-},v^n_{t^n_i-})]{}^t\delta x^n_i \delta x^n_i]
    + \sum_{i=0}^{k(n)-1} r^n_i
\ea $[\nabla^2_x F_{t^n_{i}}(x^n_{t^n_i-},v^n_{t^n_i-})]1_{t \in
]t^n_i,t^n_{i+1}]}$ is bounded by $C$, and converges to $\nabla^2_x
F_t(x_t,v_{t-})$ by left-continuity of $\nabla^2_x F$, and the paths
of both are left-continuous by lemma \ref{cadlag.prop}. Since $x$
and the subdivision $(\pi_n)$ are as in definition
\ref{finitequadraticvariation.def}, lemma
\ref{vagueconvergence2.lemma} in appendix \ref{measuretheory.sec}
applies and gives as limit: \ba \int_{0}^{T} \frac{1}{2} {\rm
tr}[{}^t\nabla^2_xF_t(x_{u},v_{u-})] d[x](u)] = \int_{0}^{T}
\frac{1}{2} {\rm tr}[{}^t\nabla^2_xF_t(x_{u},v_{u})] d[x](u)]\ea
since $\nabla^2_xF$ is predictable in the second variable i.e.
verifies \eqref{predictable.eq}. Using the same lemma, since
$|r^n_i|$ is bounded by $\epsilon^n_i |\delta x^n_i|^2$ where
$\epsilon^n_i$ converges to 0 and is bounded by $2C$,
$\sum_{i=i^n(s)+1}^{i^n(t)-1} r^n_i$ converges to 0.

\end{itemize}

Since all other terms converge, the limit:
\ba \lim_n \sum_{i=0}^{k(n)-1} \nabla_xF_{t^n_{i}}(x^n_{t^n_i-},v^n_{t^n_i-}) (x(t^n_{i+1})-x(t^n_i)) \ea
exists, and the result is established.

\end{proof}

\section{Change of variable formula for functionals of a cadlag path}

We will now extend the previous result to functionals of cadlag
paths. The following definition is a taken from F\"ollmer
\cite{follmer79}:
\begin{definition} \label{finitequadraticvariationjumps.def}
Let $\pi_n=(t^n_0,\ldots,t^n_{k(n)})$, where $0=t^n_0 \leq t^n_1
\leq \ldots \leq t^n_{k(n)}=T$ be a sequence of subdivisions of
$[0,T]$ with step decreasing to 0 as $n \rightarrow \infty$. $f \in
D([0,T],\mathbb{R})$ is said to have finite quadratic variation
along $(\pi_n)$ if the sequence of discrete measures: \ba
\xi^n=\sum_{i=0}^{k(n)-1} (f(t^n_{i+1})-f(t^n_i))^2 \delta_{t^n_i}
\ea where $\delta_t$ is the Dirac measure at $t$, converge vaguely
to a Radon measure $\xi$ on $[0,T]$ such that \ba [f](t) =\xi([0,t])
= [f]^c(t) + \sum_{0 < s \leq t} (\Delta f(s))^2 \ea where $[f]^c$
is the continuous part of $[f]$. $[f]$ is called quadratic variation
of $f$ along the sequence $(\pi_n)$. $x \in \mathcal{U}_T$ is said
to have finite quadratic variation along the sequence $(\pi_n)$ if
the functions $x_i,1 \leq i \leq d$ and $x_i+x_j, 1 \leq i < j \leq
d$ do. The quadratic variation of $x$ along $(\pi_n)$ is the
$S^+_d$-valued function $x$ defined by: \ba [x]_{ii}=[x_i],\qquad
[x]_{ij}=\frac{1}{2} ([x_i+x_j]-[x_i]-[x_j]),\quad i \neq j \ea

\end{definition}

\begin{theorem}[Change of variable formula for functionals of discontinuous paths]\label{itojumpsnoprobability.theorem}
Let $(x,v) \in \mathcal{U}_T \times \mathcal{S}_T$ where $x$ has
finite quadratic variation along $(\pi_n)$ and \ba
\label{exhaustion.eq} \sup_{t \in [0,T]-\pi_n}
|x(t)-x(t-)|+|v(t)-v(t-)| \rightarrow 0\ea Denote
\ba x^n(t)=\sum_{i=0}^{k(n)-1} x(t_{i+1}-)1_{[t_i,t_{i+1})}(t)+x(T)1_{\{T\}}(t)  \nonumber \\
v^n(t)=\sum_{i=0}^{k(n)-1}
v(t_{i})1_{[t_i,t_{i+1})}(t)+v(T)1_{\{T\}}(t), \qquad
h^n_i=t^n_{i+1}-t^n_i \label{piecewiseconstant.eq}\ea Then for any
non-anticipative functional $F \in \mathbb{C}^{1,2}$ satisfying the
following assumptions:
\begin{enumerate}
\item $F$ is predictable in the second variable in the sense of \eqref{predictable.eq}
\item $\nabla^2_x F$ and $\mathcal{D}F$ have the local boundedness  property \eqref{localbounded.eq}
\item $F,\nabla_xF,\nabla^2_xF \in \mathbb{F}^{\infty}_l$
\item $\nabla_xF$ has the horizontal local Lipschitz property \eqref{localLipschitz.eq} (Definition
\ref{timelipschitz.def})
\end{enumerate}
the following  limit exists \ba \int_{]0,T]}
\nabla_xF_t(x_{t-},v_{t-}) d^\pi x\
:=\quad\mathop{\lim}_{n\to\infty} \sum_{i=0}^{k(n)-1}
\nabla_xF_{t^n_{i}}(x^{n,\Delta x(t^n_i)}_{t^n_i-},v^n_{t^n_i-})
(x(t^n_{i+1})-x(t^n_i)) \label{follmerintegraljumps.eq} \ea and
\begin{eqnarray}
F_T(x_T,v_T)-F_0(x_0,v_0) = \int_{]0,T]}
\mathcal{D}_tF_t(x_{u-},v_{u-}) du+ \int_{]0,T]} \frac{1}{2} {\rm
tr}\left({}^t\nabla^2_xF_t(x_{u-},v_{u-})d[x]^c(u)\right)\nonumber\\
+ \int_{]0,T]} \nabla_xF_t(x_{t-},v_{t-}) d^\pi x + \sum_{u \in
]0,T]}
[F_u(x_u,v_u)-F_u(x_{u-},v_{u-})-\nabla_xF_u(x_{u-},v_{u-}).\Delta
x(u)] \label{functional.ito.jumps.eq}
\end{eqnarray}
\end{theorem}

\begin{remark} Condition \eqref{exhaustion.eq} simply means that the subdivision asymptotically contains all discontinuity points of $(x,v)$.
Since a cadlag function  has at most a countable  set of
discontinuities, this can always be achieved by adding e.g. the
discontinuity points $\{t\in [0,T], \max(|\Delta x(t)|,|\Delta
v(t)|)\geq 1/n\}$ to $\pi_n$.\end{remark}

\begin{proof}
Denote $\delta x^n_i=x(t^n_{i+1})-x(t^n_i)$. Denote also
\ba \eta_n=\sup\{|v(u)-v(t^n_i)|+|x(u)-x(t^n_i)|+|t^n_{i+1}-t^n_i|,0 \leq i \leq k(n)-1,u\in[t^n_{i},t^n_{i+1})\} \ea
and note that this quantity converges to 0 as $n \rightarrow \infty$, thanks to lemma \ref{uniform.lemma}. We assume $n$ sufficiently large so that for some $C > 0$ such that, for any $t < T$, for any $(x',v') \in \mathcal{U}_t \times \mathcal{S}_t$, $d_\infty((x_t,v_t),(x',v')) < \eta_n \Rightarrow |\mathcal{D}_tF_t(x',v')| \leq C,|\nabla^2_xF_t(x',v')| \leq C$ , using the local boundedness property satisfied by these derivatives.

For $\epsilon > 0$, we separate the jump times of $x$ in two sets: a finite set $C_1(\epsilon)$ and a set $C_2(\epsilon)$ such that $\sum_{s \in C_2(\epsilon)} |\Delta x_s|^2 \ < \epsilon^2$. We also separate the indices $0 \leq i \leq k(n)-1$ in two sets: a set $I^n_1(\epsilon)$ such that $(t_i,t_{i+1}]$ contains at least a time in $C_1(\epsilon)$, and its complementary $I^n_2(\epsilon)$.
Denoting $K=\overline{\{x(u), s \leq u \leq t \}}$ which is a compact subset of $U$, and $U^{c}=\mathbb{R^d}-U$, one may choose $\epsilon$ sufficiently small and $n$ sufficiently large so that $d(K,U^{c}) > \epsilon+\eta_n$.\\
Denote $i^n(t)$ the index such that $t \in [t^n_i,t^n_{i+1})$.
Property \eqref{exhaustion.eq} implies that for $n$ sufficiently
$C_1(\epsilon)\subset \{t^n_{i+1},i=1..k(n)\}$ so
 \ba \sum_{0 \leq i \leq k(n)-1,i \in
I^n_1(\epsilon)} F_{t^n_{i+1}}(x^{n,\Delta
x(t^n_{i+1})}_{t^n_{i+1}-},v^n_{t^n_{i+1}-})-F_{t^n_i}(x^{n,\Delta
x(t^i_n)}_{t^n_i-},v^n_{t^n_i-})  \rightarrow \sum_{u \in ]0,T] \cup
C_1(\epsilon)} F_u(x_u,v_u)-F_u(x_{u-},v_{u-}) \ea
as $n \rightarrow \infty$, by left-continuity of $F$.\\
Let us now consider, for $i \in I^n_2(\epsilon), i \leq k(n)-1$, the
decomposition:
\ba F_{t^n_{i+1}}(x^{n,\Delta x(t^n_{i+1})}_{t^n_{i+1}-},v^n_{t^n_{i+1}-})-F_{t^n_i}(x^{n,\Delta x(t^n_i)}_{t^n_i-},v^n_{t^n_i-})&=& F_{t^n_{i+1}}(x^{n,\Delta x(t^n_{i+1})}_{t^n_{i+1}-},v^n_{t^n_{i+1}-})-F_{t^n_{i+1}}(x^n_{t^n_{i+1}-},v^n_{t^n_{i+1}-}) \nonumber \\
&+&F_{t^n_{i+1}}(x^n_{t^n_{i+1}-},v^n_{t^n_i,h^n_i})-F_{t^n_i}(x^n_{t^n_i},v^n_{t^n_i}) \nonumber \\
&+&F_{t^n_i}(x^n_{t^n_{i}},v^n_{t^n_i-})-F_{t^n_{i}}(x^{n,\Delta x(t^n_i)}_{t^n_i-},v^n_{t^n_i-})
 \label{62.eq}\ea
where we have used the property \eqref{predictable.eq} to obtain
$F_{t^n_i}(x^n_{t^n_i},v^n_{t^n_i})=F_{t^n_i}(x^n_{t^n_{i}},v^n_{t^n_i-})$.
The second line in the right-hand side can be written $
\psi(h^n_i)-\psi(0) $ where: \ba
\psi(u)=F_{t^n_i+u}(x^n_{t^n_i,u},v^n_{t^n_i,u}) \ea Since $F \in
\mathbb{C}^{1,2}([0,T])$, $\psi$ is right-differentiable, and
moreover by lemma \ref{cadlag.prop}, $\psi$ is continuous, so: \ba
F_{t^n_{i+1}}(x^n_{t^n_i,h^n_i},v^n_{t^n_i,h^n_i})-F_{t^n_i}(x^n_{t^n_i},v^n_{t^n_i})=
\int_{0}^{t^n_{i+1}-t^n_i}
\mathcal{D}_{t^n_i+u}F(x^n_{t^n_i,u},v^n_{t^n_i,u}) du \ea

The third line can be written $\phi(x(t^n_{i+1}-)-x(t^n_i))-\phi(0)$, where:
\ba \phi(u)=F_{t^n_{i}}(x^{n,\Delta x(t^n_i)+u}_{t^n_i-},v^n_{t^n_i-}) \ea
Since $F \in \mathbb{C}^{1,2}([0,T])$, $\phi$ is well-defined and $C^2$ on the convex set $B(x(t^n_i),\eta_n+\epsilon) \subset U$, with:
\ba \phi'(u)=\nabla_x F_{t^n_{i}}(x^{n,\Delta x(t^n_i)+u}_{t^n_i-},v^n_{t^n_i-})  
\phi''(u)=\nabla^2_x F_{t^n_{i}}(x^{n,\Delta x(t^n_i)+u}_{t^n_i-},v^n_{t^n_i-}) \ea
So a second order Taylor expansion of $\phi$ at $u=0$ yields:
\ba F_{t^n_{i}}(x^n_{t^n_{i}},v^n_{t^n_i-})-F_{t^n_{i}}(x^{n,\Delta x(t^n_i)}_{t^n_i-},v^n_{t^n_i-})=\nabla_x F_{t^n_{i}}(x^{n,\Delta x(t^n_i)}_{t^n_i-},v^n_{t^n_i-})(x(t^n_{i+1}-)-x(t^n_i)) \nonumber \\
+ \frac{1}{2} tr[\nabla^2_x F_{t^n_{i}}(x^{n,\Delta
x(t^n_i)}_{t^n_i-},v^n_{t^n_i-})]{}^t(x(t^n_{i+1}-)-x(t^n_i))
(x(t^n_{i+1}-)-x(t^n_i))] +r^n_i \ea where $r^n_{i,1}$ is bounded by
\ba K |(x(t^n_{i+1}-)-x(t^n_i))|^2 \sup_{x \in
B(x(t^n_i),\eta_n+\epsilon)} |\nabla^2_x
F_{t^n_{i}}(x^{n,x-x(t^n_i)}_{t^n_i-},v^n_{t^n_i-})-\nabla^2_x
F_{t^n_{i}}(x^n_{t^n_i-},v^n_{t^n_i-})| \ea Similarly, the first
line can be written
$\phi(x(t^n_{i+1})-x(t_i^n))-\phi(x(t^n_{i+1}-)-x(t^n_i))$ where
$\phi(u)=F_{t^n_{i+1}}(x^{n,\Delta
x(t^n_i)+u}_{t^n_{i}-,h^n_i},v^n_{t^n_{i-}})$. So, a second order
Taylor expansion of $\phi$ at $u=0$ yields:
\ba F_{t^n_{i+1}}(x^{n,\Delta x(t^n_{i+1})}_{t^n_{i+1}-},v^n_{t^n_{i+1}-})-F_{t^n_{i+1}}(x^n_{t^n_{i+1}-},v^n_{t^n_{i+1}-})=\nabla_x F_{t^n_{i+1}}(x^{n,\Delta x(t^n_i)}_{t^n_i-,h^n_i},v^n_{t^n_i,h^n_i})\Delta x(t^n_{i+1}) \nonumber \\
+ \frac{1}{2} tr[\nabla^2_x F_{t^n_{i+1}}(x^{n,\Delta
x(t^n_i)}_{t^n_i-,h^n_i},v^n_{t^n_i,h^n_i})]{}^t\Delta x(t^n_{i+1})
\Delta x(t^n_{i+1})+r^n_{i,2} \ea where $r^n_{i,2}$ is bounded by
\ba K |\Delta x(t^n_{i+1})|^2 \sup_{x \in
B(x(t^n_i),\eta_n+\epsilon)} |\nabla^2_x
F_{t^n_{i}}(x^{n,x-x(t^n_i-)}_{t^n_i-},v^n_{t^n_i-})-\nabla^2_x
F_{t^n_{i}}(x^{n,\Delta x(t^n_i)}_{t^n_i-},v^n_{t^n_i-})| \ea Using
the horizontal local Lipschitz property \eqref{localLipschitz.eq}
for $\nabla_xF$, for $n$ sufficiently large: \ba |\nabla_x
F_{t^n_{i+1}}(x^{n,\Delta
x(t^n_i)}_{t^n_i-,h^n_i},v^n_{t^n_i,h^n_i})-\nabla_x
F_{t^n_{i}}(x^{n,\Delta
x(t^n_i)}_{t^n_i-},v^n_{t^n_i-})|<C(t^n_{i+1}-t^n_i) \ea On other
hand, since $\nabla^2_xF$ is bounded by $C$ on all paths considered:
\ba |\quad{\rm tr}\left(\nabla^2_x F_{t^n_{i}}(x^{n,\Delta x(t^n_i)}_{t^n_i-},v^n_{t^n_i-}){}^t(x(t^n_{i+1}-)-x(t^n_i)) (x(t^n_{i+1}-)-x(t^n_i))\right) \nonumber \\
+ {\rm tr}\left(\nabla^2_x F_{t^n_{i+1}}(x^{n,\Delta x(t^n_i)}_{t^n_i-,h^n_i},v^n_{t^n_i,h^n_i})]{}^t\Delta x(t^n_{i+1}) \Delta x(t^n_{i+1})\right) \nonumber \\
-{\rm tr}\left(\nabla^2_x F_{t^n_{i}}(x^{n,\Delta
x(t^n_i)}_{t^n_i-},v^n_{t^n_i-})]{}^t\delta x^n_i \delta
x^n_i\right) | < 2C| \Delta x(t^n_{i+1})|^2 \ea Hence, we have shown
that:
\ba F_{t^n_{i+1}}(x^{n,\Delta x(t^n_{i+1})}_{t^n_{i+1}-},v^n_{t^n_{i+1}-})-F_{t^n_{i+1}}(x^n_{t^n_{i+1}-},v^n_{t^n_{i+1}-})+F_{t^n_{i}}(x^n_{t^n_{i}},v^n_{t^n_i-})-F_{t^n_{i}}(x^{n,\Delta x(t^n_i)}_{t^n_i-},v^n_{t^n_i-}) = \nonumber \\
\nabla_x F_{t^n_{i}}(x^{n,\Delta
x(t^n_i)}_{t^n_i-},v^n_{t^n_i-})\delta x^n_i + \frac{1}{2}
tr[\nabla^2_x F_{t^n_{i}}(x^{n,\Delta
x(t^n_i)}_{t^n_i-},v^n_{t^n_i-})]{}^t \delta x^n_i \delta x^n_i]
+r^n_i+q^n_i \nonumber \ea where $r^n_i$ is bounded by: \ba  4K
|\delta x^n_i|^2 \sup_{x \in B(x(t^n_i),\eta_n+\epsilon)}
|\nabla^2_x
F_{t^n_{i}}(x^{n,x-x(t^n_i-)}_{t^n_i-},v^n_{t^n_i-})-\nabla^2_x
F_{t^n_{i}}(x^{n,\Delta x(t^n_i)}_{t^n_i-},v^n_{t^n_i-})| \ea and
$q^n_i$ is bounded by: \ba C'( h^n_i|\Delta x(t^n_i)|+|\Delta
x(t^n_i)|^2 )\ea

Denote $i^n(t)$ the index such that $t \in
[t^n_{i^n(t)},t^n_{i^n(t)+1}[$. Summing all the terms  above for
  $i \in C_2(\epsilon)\cap \{0,1,..k(n)-1\}$:
\begin{itemize}
\item The left-hand side of \eqref{62.eq}  yields \ba F_{T}(x^n_{T},v^n_{T})-F_{0}(x_0,v_0)
    -\sum_{0\leq i \leq k(n)-1,i \in I^n_1(\epsilon)} F_{t^n_{i+1}}(x^n_{t^n_{i+1}},v^n_{t^n_{i+1}})-F_{t^n_i}(x^n_{t^n_i},v^n_{t^n_i}) \ea
    which converges to
    \ba F_{T}(x_T,v_T)-F_0(x_0,v_0)-\sum_{u \in ]0,T] \cup C_1(\epsilon)} F_u(x_u,v_u)-F_u(x_{u-},v_{u-}) \ea

\item The sum of the first and third lines of \eqref{62.eq} the right-hand side can be written:
\ba \sum_{0 \leq i \leq k(n)-1, i \in I^n_2(\epsilon)} \nabla_x F_{t^n_{i}}(x^{n,\Delta x(t^n_i)}_{t^n_i-},v^n_{t^n_i-})\delta x^n_i  \nonumber \\
    + \sum_{0 \leq i \leq k(n)-1, i \in I^n_2(\epsilon)} \frac{1}{2} {\rm tr}\left(\nabla^2_x F_{t^n_{i}}(x^{n,\Delta x(t^n_i)}_{t^n_i-},v^n_{t^n_i-}){}^t \delta x^n_i \delta x^n_i\right) \nonumber \\
    + \sum_{0 \leq i \leq k(n)-1, i \in I^n_2(\epsilon)} r^n_i+q^n_i \label{63.eq}
\ea Consider the measures $\mu^n_{ij}=\xi^n_{ij}-\sum_{0 < s \leq T,
s \in C_2(\epsilon)} (\Delta f_{ij}(s))^2 \delta_{s}$, where
$f_{ii}=x_i,1 \leq j \leq d$ and $f_{ij}=x_i+x_j, 1 \leq i < j \leq
d$ and $\xi^n_{ij}$ is defined in Definition
\ref{finitequadraticvariationjumps.def}. The second line of
\eqref{63.eq} can be  decomposed as: \ba A_n + \frac{1}{2}\sum_{0 <
u \leq T, u \in C_2(\epsilon)} {\rm tr}\left( \nabla^2_x
F_{t^n_{i}}(x^{n,\Delta x(t^n_i)}_{t^n_i-},v^n_{t^n_i-}){}^t \Delta
x(u) \Delta x(u)\right) \label{61.eq}\ea where
$$ A_n = {\rm tr}\ \int_{]0,T]} \mu^n(dt)\quad \sum_{0 \leq i \leq k(n)-1, i \in I^n_2(\epsilon)}\nabla^2_x
F_{t^n_{{i^n(t)}}}(x^{n,\Delta
x(t^n_{i^n(t)})}_{t^n_{i^n(t)}-},v^n_{t^n_{i^n(t)}-})\quad 1_{t \in
(t^n_i,t^n_{i+1}]} $$ where $\mu^n$ denotes the matrix-valued
measure with components  $\mu^n_{ij}$ defined above.  $\mu^n_{ij}$
 converges vaguely to
the atomless measure $[f_{ij}]^c$. Since $\sum_{0 \leq i \leq
k(n)-1, i \in I^n_2(\epsilon)} \nabla^2_x
F_{t^n_{{i^n(t)}}}(x^{n,\Delta
x(t^n_{i^n(t)})}_{t^n_{i^n(t)}-},v^n_{t^n_{i^n(t)}-})\ 1_{t \in
(t^n_i,t^n_{i+1}]}$ is bounded by $C$ and converges to $\nabla^2_x
F_t(x_{t-},v_{t-})1_{t\notin C_1(\epsilon)}$ by left-continuity of
$\nabla^2_x F$,  applying Lemma \ref{vagueconvergence2.lemma}  to
$A_n$ and yields  that $A_n$ converges to: \ba \int_{]0,T]}
\frac{1}{2} {\rm
tr}\left({}^t\nabla^2_xF_t(x_{u-},v_{u-})d[x]^c(u)\right) \ea The
second term  in \eqref{61.eq} has the $\limsup$ of its absolute
value bounded by $C \epsilon^2$. Using the same argument, since
$|r^n_i|$ is bounded by $s^n_i |\delta x^n_i|^2$ for some $s^n_i$
which converges to 0 and is bounded by some constant,
$\sum_{i=0}^{k(n)-1} |r^n_i|$ has its $\limsup$ bounded by $2C
\epsilon^2$; similarly, the $\limsup$ of $\sum_{i=0}^{k(n)-1}
|q^n_i|$ is bounded by $C'(T\epsilon + \epsilon^2)$.

The term in the first line of \eqref{63.eq} can be written:
\ba \sum_{i=0}^{k(n)-1} \nabla_x F_{t^n_{i}}(x^{n,\Delta x(t^n_i)}_{t^n_i-},v^n_{t^n_i-}) (x(t^n_{i+1})-x(t^n_i)) \nonumber \\
- \sum_{0 \leq i \leq k(n)-1, i \in I^n_1 (\epsilon)} \nabla_x F_{t^n_{i}}(x^{n,\Delta x(t^n_i)}_{t^n_i-},v^n_{t^n_i-}) (x_{t^n_{i+1}}-x_{t^n_i}) \ea
where the second term converges to $\sum_{0 < u \leq T, u \in C_1(\epsilon)} \nabla_xF_u (x_{u-},v_{u-}) \Delta x(u)$.

\item The second line of \eqref{62.eq}:
\ba \int_{0}^{T}
\mathcal{D}_tF_{u}(x^n_{t^n_{i^n(u)},u-t^n_{i^n(u)}},v^n_{t^n_{i^n(u)},u-t^n_{i^n(u)}})1_{i^n(u)
\in I^n_2(\epsilon)} du \ea where the integrand converges to
$\mathcal{D}_tF_{u}(x_{u-},v_{u-})1_{u \notin C_1(\epsilon)}$ and is
bounded by $C$, hence by dominated convergence this term converges
to: \ba \int_{0}^{T} \mathcal{D}_tF_t(x_{u-},v_{u-}) du \ea

Summing up, we have established that the difference between the $\limsup$ and the $\liminf$ of:
\ba \sum_{i=0}^{k(n)-1} \nabla_x F_{t^n_{i}}(x^{n,\Delta x(t^n_i)}_{t^n_i-},v^n_{t^n_i-}) (x(t^n_{i+1})-x(t^n_i))1_{s < t^n_i \leq t} \ea
is bounded by $C'' (\epsilon^2+T \epsilon)$. Since this is true for any $\epsilon$, this term has a limit.

Let us now write the equality we obtained for a fixed $\epsilon$:
\ba F_T(x_{T},v_{tT})-F_0(x_0,v_0) = \int_{]0,T]}
\mathcal{D}_tF_t(x_{u-},v_{u-}) du
+ \int_{]0,T]} \frac{1}{2} {\rm tr}[{}^t\nabla^2_xF_t(x_{u-},v_{u-})d[x]^c(u)] \nonumber \\
+ \lim_n \sum_{i=0}^{k(n)-1}\nabla_x F_{t^n_{i}}(x^{n,\Delta
x(t^n_i)}_{t^n_i-},v^n_{t^n_i-}) (x(t^n_{i+1})-x(t^n_i))
\nonumber\\+ \sum_{u \in ]0,T] \cup C_1(\epsilon)}
[F_u(x_u,v_u)-F_u(x_{u-},v_{u-})-\nabla_xF_u(x_{u-},v_{u-}) \Delta
x(u)] + \alpha(\epsilon) \nonumber \ea where $\alpha(\epsilon) \leq
C'' (\epsilon^2+T \epsilon)$. The only point left to show   is that:
\ba \sum_{u \in ]0,T] \cup C_1(\epsilon)}
[F_u(x_u,v_u)-F_u(x_{u-},v_{u-})-\nabla_xF_u(x_{u-},v_{u-}) \Delta
x(u)] \ea converges to: \ba \sum_{u \in ]0,T]}
[F_u(x_u,v_u)-F_u(x_{u-},v_{u-})-\nabla_xF_u(x_{u-},v_{u-}) \Delta
x(u)] \ea
which is to say that the sum above is absolutely convergent.\\
We can first choose $d(K,U^{c}) > \eta > 0$ such that: \ba \forall u
\in [0,T], \quad \forall (x',v') \in \mathcal{U}_u \times
\mathcal{S}_u, d_{\infty}((x_t,v_t),(x',v')) \leq \eta \Rightarrow
|\nabla^2_x F_u(x(u),v(u))| < C
 \ea
The jumps of $x$ of magnitude greater than $\eta$ are in finite
number. Then, if $u$ is a jump time of $x$ of magnitude less than
$\eta$, then $x(u-)+h\Delta x(u) \in U$ for $h \in [0,1]$, so that
we can write:
\ba F_u(x_u,v_u)-F_u(x_{u-},v_{u-})-\nabla_xF_u(x_{u-},v_{u-}) \Delta x(u) = \nonumber \\
\int_0^1 (1-v)[{}^t \nabla^2_xF_u(x^{h \Delta x(u)}_{u-},v_{u-}) {}^t \Delta x(u) \Delta x(u)] \nonumber \
\leq \frac{1}{2}C |\Delta x(u)|^2 \ea

Hence, the theorem is established.

\end{itemize}

\end{proof}

\begin{remark} \label{rightcontinuousito.remark}
If the vertical derivatives are right-continuous instead of left-continuous, and $\nabla_xF$ not necessarily locally Lipschitz in time, define:
\ba x^n(t)=\sum_{i=0}^{k(n)-1} x(t_{i})1_{[t_i,t_{i+1})}(t)+x(T)1_{\{T\}}(t)  \nonumber \\
v^n(t)=\sum_{i=0}^{k(n)-1} v(t_{i})1_{[t_i,t_{i+1})}(t)+v(T)1_{\{T\}}(t) \qquad h^n_i=t^n_{i+1}-t^n_i \ea
Following the same argument than in the proof with the decomposition:
\ba F_{t^n_{i+1}}(x^n_{t^n_{i+1}},v_{t^n_{i+1}})-F_{t^n_i}(x^n_{t^n_i},v_{t^n_i})&=&F_{t^n_{i+1}}(x^n_{t^n_{i+1}},v_{t^n_{i+1}})-F_{t^n_{i+1}}(x^n_{t^n_{i+1}},v_{t^n_i,h^n_i}) \nonumber \\
&+&F_{t^n_{i+1}}(x^n_{t^n_{i+1}},v_{t^n_i,h^n_i})-F_{t^n_{i+1}}(x^n_{t^n_i,h^n_i},v_{t^n_i,h^n_i}) \nonumber \\
&+&F_{t^n_{i+1}}(x^n_{t^n_i,h^n_i},v_{t^n_i,h^n_i})-F_{t^n_i}(x^n_{t^n_i},v_{t^n_i}) \ea
leads to the formula with the following Riemann sum:
\ba \lim_n \sum_{i=0}^{k(n)-1} \nabla_xF_{t^n_{i+1}}(x^n_{t^n_i,h^i},v^n_{t^n_i,h^n_i}) (x(t^n_{i+1})-x(t^n_i)) \ea
\end{remark}


\section{Functionals of Dirichlet processes} \label{dirichlet.sec}

A   Dirichlet process \cite{follmer81,coquet03}, or finite energy
process, on a filtered probability space
$(\Omega,\mathcal{B},(\mathcal{B}_t),\mathbb{P})$ is an adapted
cadlag process  that can be represented as the sum of a
semimartingale and an adapted continuous process with zero quadratic
variation along dyadic subdivisions.

 For continuous Dirichlet
processes, a pathwise It\^o calculus was introduced by   H.
F{\"o}llmer in \cite{follmer79,follmer81,lyons94}. Coquet, M\'emin
and Slominski \cite{coquet03} extended these
 results to discontinuous Dirichlet processes \cite{stricker88}.
Using Theorem \ref{itojumpsnoprobability.theorem} we can
 extend these results to functionals of Dirichlet processes;
 this yields in particular  a pathwise construction of stochastic integrals for functionals of a Dirichlet process.

Let $Y(t)=X(t)+B(t)$ be a $U$-valued Dirichlet process defined as
the sum of a semimartingale $ X  $ on some filtered probability
space $(\Omega,{\cal B},{\cal B}_t,\mathbb{P})$ and $B$ an adapted
continuous process $B$ with zero quadratic variation along the
dyadic subdivision. We denote by $[X]$ the quadratic variation
process associated to $X$, $[X]^c$ the continuous part of $[X]$, and
$\mu(dt\ dz)$ the integer-valued random  measure describing the
jumps of $X$ (see \cite{jacodshiryaev} for definitions).

Let $A$ be an adapted  process with $S$-valued   cadlag paths. Note
that $A$ need not be  a semimartingale.

We call $\Pi_n = \{ 0=t^n_0 < t^n_1 < \ldots < t^n_{k(n)}=T\}$  a
{\it random subdivision} if the $t^n_i$ are stopping times with
respect to $({\cal B}_t)_{t\in[0,T]}$.
\begin{proposition}[Change of variable formula for Dirichlet processes] \label{dirichlet.prop}

Let  $\Pi_n = \{ 0= t^n_0 < t^n_1 < \ldots < t^n_{k(n)}=T\}$ be any
sequence of
random subdivisions of $[0,T]$   such that \\
(i) $X$ has finite quadratic variation along $\Pi_n$ and $B$ has
zero quadratic variation along $\Pi_n$ almost-surely,$$ (ii)\qquad
\mathop{\sup}_{t \in [0,T]-\Pi_n} |Y(t)-Y(t-)|+|A(t)-A(t-)|
\mathop{\rightarrow}^{n\to\infty} 0\quad \mathbb{P}-a.s.$$
 Then there exists $\Omega_1\subset \Omega$ with $\mathbb{P}(\Omega_1)=1$ such that for any
non-anticipative functional $F\in \mathbb{C}^{1,2}$ satisfying
\begin{enumerate}
\item $F$ is predictable in the second variable in the sense of
\eqref{predictable.eq}
\item $\nabla^2_x F$ and $\mathcal{D}F$ satisfy the local boundedness property  \eqref{localbounded.eq}
\item $F,\nabla_xF,\nabla^2_xF \in \mathbb{F}^{\infty}_l$
\item $\nabla_xF$ has the horizontal  local Lipschitz property
\eqref{localLipschitz.eq},
\end{enumerate}
the following equality holds on $\Omega_1$ for all $t \leq T$: \ba
F_t(Y_t,A_t)-F_0(Y_0,A_0) =
\int_{]0,t]} \mathcal{D}_uF(Y_{u-},A_{u-}) du+ \int_{]0,t]} \frac{1}{2} {\rm tr}[{}^t\nabla^2_xF_u(Y_{u-},A_{u-}) d[X]^c(u)] \nonumber \\
+ \int_{]0,t]} \int_{R^d}
[F_u(Y^{z}_{u-},{A_u-})-F_u(Y_{u-},A_{u-})-z\nabla_xF_u(Y_{u-},A_{u-})]\mu(du,dz)\nonumber\\
+ \int_{]0,t]} {\nabla}_xF_u(Y_{u-},A_{u-}).dY(u)  \quad a.s.
\label{dirichletito.eq}\ea where the last term is the F\"ollmer
integral \eqref{follmerintegraljumps.eq} along the subdivision
$\Pi_n$, defined for $\omega\in\Omega_1$ by: \ba \quad\int_{]0,t]}
{\nabla}_xF_u(Y_{u-},A_{u-}).dY(u) := \lim_n \sum_{i=0}^{k(n)-1}
\nabla_xF_{t^n_{i}}(Y^{n,\Delta Y(t^n_i)}_{t^n_i-},A^n_{t^n_i-})
(Y(t^n_{i+1})-Y(t^n_i))1_{]0,t]} \label{follmerdirichlet.eq}\ea
where $(Y^n,A^n)$ are the piecewise constant approximations along
$\Pi_n$, defined as in \eqref{piecewiseconstant.eq}.

Moreover, the F\"ollmer integral with respect to any other random
subdivision verifying (i)--(ii),  is almost-surely equal to
\eqref{follmerdirichlet.eq}.
\end{proposition}
\begin{remark} Note that the convergence of
\eqref{follmerdirichlet.eq} holds over a set $\Omega_1$ which may be
chosen independently of the choice of $F\in \mathbb{C}^{1,2}$.
\end{remark}
\begin{proof}
Let $(\Pi_n)$ be a sequence of random subdivisions verifying
(i)--(ii). Then there exists a set  $\Omega_1$ with
$\mathbb{P}(\Omega_1)=1$ such that for  $\omega\in\Omega_1$ $(X,A)$
is a  cadlag function and (i)-(ii) hold pathwise. Applying Theorem
\ref{itojumpsnoprobability.theorem} to $(Y(.,\omega),A(.,\omega))$
along the subdivision $\Pi_n(\omega)$ shows that
\eqref{dirichletito.eq} holds on $\Omega_1$.

To show independence of the limit in \eqref{follmerdirichlet.eq}
from the chosen subdivision, we note that if $\Pi^2_n$ another
sequence of random subdivisions satisfies (i)--(ii), there exists
$\Omega_2\subset \Omega$ with $\mathbb{P}(\Omega_2)=1$ such that one
can apply Theorem \ref{itojumpsnoprobability.theorem} pathwise for
$\omega\in \Omega_2$. So we have
$$ \int_{]0,t]}
{\nabla}_xF_u(Y_{u-},A_{u-}).d^{\Pi^2}Y(u)=\int_{]0,t]}
{\nabla}_xF_u(Y_{u-},A_{u-}).d^{\Pi}Y(u)$$ on $\Omega_1\cap
\Omega_2$. Since $\mathbb{P}(\Omega_1\cap \Omega_2)=1$ we obtain the
result.
\end{proof}

\section{Functionals of semimartingales}
\label{semimartingale.sec} Proposition \ref{dirichlet.prop} holds of
course when $X$ is a semimartingale. We will now show that in this
case,
 under an additional assumption,  the pathwise  integral  coincides    almost-surely with the stochastic integral $\int Y dX$.
\subsection{Cadlag semimartingales}
Let $X$ be a cadlag semimartingale and $A$ an adapted cadlag process
on $(\Omega,{\cal B},{\cal B}_t,\mathbb{P})$. We use the notations
$[X]$ , $[X]^c$, $\mu(dt\ dz)$ defined in Section
\ref{dirichlet.sec}.

 Theorem \ref{itojumpsnoprobability.theorem}  yields an
It\^o formula for functionals of  $X$: under the additional
assumption $\nabla_xF \in \mathbb{B}$,  the pathwise F\"ollmer
integral coincides with the  stochastic integral.

\begin{proposition}[Functional It\^o formula for a semimartingale] \label{itojumps.prop}
Let $F \in \mathbb{C}^{1,2}$ be a non-anticipative functional
satisfying
\begin{enumerate}
\item $F$ is predictable in the second variable, i.e. verifies
\eqref{predictable.eq},
\item $\nabla_xF,\ \nabla^2_x F,\  \mathcal{D}F\in \mathbb{B}$,
\item $F,\nabla_xF,\nabla^2_xF \in \mathbb{F}^{\infty}_l$,
\item $\nabla_xF$ has the horizontal local Lipschitz property
\ref{localLipschitz.eq}.
\end{enumerate}
Then: \begin{align} F_t(X_t,A_t)-F_0(X_0,A_0) =&\nonumber\\
\int_{]0,t]} \mathcal{D}_uF(X_{u-},A_{u-}) du+& \int_{]0,t]}
\frac{1}{2} {\rm tr}[{}^t\nabla^2_xF_u(X_{u-},A_{u-}) d[X]^c(u)]
+ \int_{]0,t]} {\nabla}_xF_u(X_{u-},A_{u-}).dX(u) \nonumber \\
+ \int_{]0,t]} \int_{R^d} &
[F_u(X^{z}_{u-},{A_u-})-F_u(X_{u-},A_{u-})-z.\nabla_xF_u(X_{u-},A_{u-})]\mu(du,dz),
\mathbb{P}\mbox{-a.s.} \end{align} where the stochastic integral is
the It\^o integral with respect to a semimartingale.

In particular,  $Y(t)=F_t(X_t,A_t)$ is a semimartingale.
\end{proposition}
\begin{remark} These results yield a non-probabilistic proof
for functional Ito formulas obtained for continuous semimartingales
 \cite{ContFournie09a,ContFournie09b,dupire09} using probabilistic
methods.
\end{remark}
\begin{proof}
Assume first that the process $X$ does not exit a compact set $K
\subset U$, and that $A$ is bounded by some constant $R > 0$. We
define the following sequence of stopping times:
\begin{eqnarray}
\tau^{n}_{0}=0 \nonumber \\
\tau^{n}_{k}=\inf\{u > \tau^{n}_{k-1} | 2^{n}u \in \mathbb {N}\
{\rm or}\quad |A(u)-A(u-)|\vee |X(u)-X(u-)| > \frac{1}{n} \}\wedge T
\end{eqnarray}
Then the coordinate processes $X_i$ and their sums $X_i+X_j$ satisfy
the property: \ba \sum_{\tau_i < s} (Z(\tau_{i})-Z(\tau_{i-1}))^2
\mathop{\rightarrow}^{\mathbb{P}}_{n \rightarrow \infty} [Z](s)
\label{bracketproba.eq} \ea in probability. There exists  a
subsequence of subdivisions such that the convergence happens almost
surely for all $s$ rational, and hence it happens almost surely for
all $s$ because both sides of \eqref{bracketproba.eq} are
right-continuous. Let $\Omega_1$ be the set on which this
convergence happens, and on which the paths of $X$ and $A$ are
$U$-valued cadlag functions. For $\omega\in \Omega_1$,  Theorem
\ref{itojumpsnoprobability.theorem} applies and yields
 \ba
F_t(X_t,A_t)-F_0(X_0,A_0) = \int_{]0,t]}
\mathcal{D}_uF(X_{u-},A_{u-}) du+ \int_{]0,t]}
\frac{1}{2} {\rm tr}[{}^t\nabla^2_xF_u(X_{u-},A_{u-}) d[X]^c(u)] \\
+ \int_{]0,t]} \int_{R^d}
[F_u(X^{z}_{u-},{A_u-})-F_u(X_{u-},A_{u-})-z.\nabla_xF_u(X_{u-},A_{u-})]\mu(du,dz)\nonumber\\
+ \mathop{\lim}_{n\to\infty}\sum_{i=0}^{k(n)-1}
\nabla_xF_{\tau^n_{i}}(X^{n,\Delta
X(\tau^n_i)}_{\tau^n_i-},A^n_{\tau^n_i-})
(X(\tau^n_{i+1})-X(\tau^n_i)) \nonumber \ea It remains to show that
the last term, which may also be written as \ba
\mathop{\lim}_{n\to\infty} \int_{]0,t]} \sum_{i=0}^{k(n)-1}
1_{]\tau^n_i,\tau^n_{i+1}]}(t)\quad
\nabla_xF_{\tau^n_{i}}(X^{n,\Delta
X(\tau^n_i)}_{\tau^n_i-},A^n_{\tau^n_i-})  .dX(t) \label{81.eq}\ea
 coincides with the (Ito)
stochastic integral of $\nabla_xF(X_{u-},A_{u-})$ with respect to
the semimartingale $X$.

First, we note that since $X,A$ are bounded and
$\nabla_xF\in\mathbb{B}$, $\nabla_xF(X_{u-},A_{u-})$ is a bounded
predictable process (by Theorem \ref{Measurability.thm}) hence its
stochastic integral $\int_0^. \nabla_xF(X_{u-},A_{u-}).dX(u)$
 is  well-defined.
  Since the integrand in \eqref{81.eq} converges almost surely to $\nabla_xF_t(X_{t-},A_{t-})$,
   and is bounded independently of $n$ by a deterministic constant $C$,
    the dominated convergence theorem for
     stochastic integrals \cite[Ch.IV Theorem32]{protter} ensures that \eqref{81.eq}
     converges in probability to $\int_{]0,t]} {\nabla}_xF_u(X_{u-},A_{u-}).dX(u)$.
     Since it converges almost-surely by proposition \ref{dirichlet.prop}, by almost-sure uniqueness of
     the limit in probability, the
      limit has to be $\int_{]0,t]} {\nabla}_xF_u(X_{u-},A_{u-}).dX(u)$.

Now we consider the general case where $X$ and $A$ may be unbounded.
Let $U^c=\mathbb{R}^d - U$ and  denote $\tau_n=\inf \{s < t|
d(X(s),U^{c}) \leq \frac{1}{n} \mbox{ or } |X(s)| \geq n \mbox{ or }
|A(s)| \geq n \}\wedge t$, which are stopping times. Applying the
previous result to the stopped processes
$(X^{\tau_n-},A^{\tau_n-})=(X(t \wedge \tau_n-),A(t \wedge
\tau_n-))$ leads to: \ba F_t(X^{\tau_n-}_t,A^{\tau_n-}_t) &=&
\int_{]0,\tau^n)} [\mathcal{D}_uF(X_u,A_u) du+ \frac{1}{2} {\rm
tr}[{}^t\nabla^2_xF_u(X_u,A_u) d[X]^c(u)] \nonumber \\
&+& \int_{]0,\tau_n)}
{\nabla}_xF_u(X_u,A_u).dX(u) \nonumber \\
&+& \int_{]0,\tau_n)} \int_{R^d} [F_u(X^{x}_{u-},{A_u-})-F_u(X_{u-},A_{u-})-z.\nabla_xF_u(X_{u-},A_{u-})]\mu(du\ dz) \nonumber \\
&+& \int_{(\tau^n,t)} {\cal D}_uF(X^{\tau_n}_u,A^{\tau_n}_u)du  \ea
Since almost surely $t\wedge\tau_n=t$ for $n$ sufficiently large,
taking the limit $n \rightarrow \infty$ yields: \ba
F_t(X_{t-},A_{t-}) &=& \int_{]0,t)} [\mathcal{D}_uF(X_u,A_u) du+
\frac{1}{2} {\rm
tr}\left({}^t\nabla^2_xF_u(X_u,A_u) d[X]^c(u)\right)  \nonumber \\
&+& \int_{]0,t)}
{\nabla}_xF_u(X_u,A_u).dX(u) \nonumber \\
&+& \int_{]0,t)} \int_{R^d}
[F_u(X^{x}_{u-},{A_u-})-F_u(X_{u-},A_{u-})-z.\nabla_xF_u(X_{u-},A_{u-})]\mu(du\
dz) \nonumber \\ \ea Adding the jump
$F_t(X_{t},A_{t})-F_t(X_{t-},A_{t-})$  to both the left-hand side
and the third line of the right-hand side, and adding
$\nabla_xF_t(X_{t-},A_{t-})\Delta X(t)$ to the second line and
subtracting it from the third, leads to the desired result.
\end{proof}

\begin{example}[Dol\'eans exponential]
Let $X$ be a scalar cadlag semimartingale, such that the continuous
part of its quadratic variation can be represented as: \ba
[X]^c(t)=\int_0^t A(s)ds \ea for some cadlag adapted process $A$.
Consider the non-anticipative functional: \ba
F_t(x_t,v_t)=e^{x(t)-\frac{1}{2}\int_0^t v(s)ds}\prod_{s \leq t}
(1+\Delta x(s))e^{- \Delta x(s)} \ea Then $F\in
\mathbb{C}^{1,\infty}$  with: \ba
\mathcal{D}_tF(x_t,v_t)=-\frac{1}{2}v(t)F_t(x_t,v_t) \ea and \ba
\nabla^k_xF_t(x_t,v_t)=F_t(x_t,v_t), k \geq 1 \ea and satisfies the
assumptions of Proposition \ref{itojumps.prop}. The process \ba
Y(t)=F_t(X_t,A_t)=e^{X(t)-\frac{1}{2}[X]^c(t)}\prod_{s \leq t}
(1+\Delta X(s))e^{- \Delta X(s)} \ea is the Dol\'eans  exponential
of the semimartingale $X$ and  Proposition \ref{itojumps.prop}
yields the well-known relation $$ Y(t)=\int_0^t Y(s-)dX(s).$$
\end{example}

\subsection{Continuous semimartingales}
In the case of a continuous semimartingale $X$ and a continuous
adapted process $A$,  an It\^o formula may also be obtained for
functionals whose  vertical derivative isright-continuous rather
than left-continuous.

\begin{proposition}[Functional It\^o formula for a continuous semimartingale] \label{itorightcontinuousderivative.prop}
Let $X$ be a continuous semimartingale with quadratic variation
process $[X]$, and $A$ a continuous adapted process, on some
filtered probability space  $(\Omega,{\cal B},{\cal
B}_t,\mathbb{P})$. Then for any non-anticipative functional $F \in
\mathbb{C}^{1,2}$
 satisfying
\begin{enumerate}
\item $F$ is predictable in the second variable, i.e. verifies
\eqref{predictable.eq},
\item $\nabla_xF,\ \nabla^2_x F,\  \mathcal{D}F\in \mathbb{B}$,
\item $F,\nabla_xF,\nabla^2_xF \in \mathbb{F}^{\infty}_l$,
\item $F \in \mathbb{F}^{\infty}_l$
\item $\nabla_xF,\nabla^2_xF \in \mathbb{F}^{\infty}_r$
\end{enumerate}
we have \ba F_t(X_t,A_t)-F_0(X_0,A_0) =
\int_0^t \mathcal{D}_uF(X_{u},A_{u}) du\nonumber\\
+ \int_0^t \frac{1}{2} {\rm tr}[{}^t\nabla^2_xF_u(X_{u},A_{u})
d[X](u)]  + \int_0^t {\nabla}_xF_u(X_{u},A_{u}).dX(u),\qquad
\mathbb{P}\mbox{-a.s.} \nonumber \ea where last term is the It\^o
stochastic integral with respect to the $X$.
\end{proposition}

\begin{proof}
Assume first that $X$ does not exit a compact set $K \subset U$ and
that $A$ is bounded by some constant $R > 0$. Let $0=t^n_0 \leq
t^n_1 \ldots \leq t^n_{k(n)}=t$ be a deterministic subdivision of
$[0,t]$. Define the approximates $(X^n,A^n)$ of $(X,A)$ as in remark
\ref{rightcontinuousito.remark}, and notice that, with the same
notations: \ba \sum_{i=0}^{k(n)-1}
\nabla_xF_{t^n_{i+1}}(X^n_{t^n_i,h^n_i},A^n_{t^n_i,h^n_i})
(X(t^n_{i+1})-X(t^n_i))=\int_{]0,t]}
\nabla_xF_{t^n_{i+1}}(X^n_{t^n_i,h^n_i},A^n_{t^n_i,h^n_i})1_{]t^n_i,t^n_{i+1}]}(t)dX(t)
\nonumber \ea
which is a well-defined stochastic integral since the integrand is predictable (left-continuous and adapted by theorem \ref{Measurability.thm}), since the times $t^n_i$ are \textit{deterministic}; this would not be the case if we had to include jumps of $X$ and/or $A$ in the subdivision as in the case of the proof of proposition \ref{itojumps.prop}. By right-continuity of $\nabla_xF$, the integrand converges to $\nabla_xF_t(X_t,A_t)$. It is moreover bounded independently of $n$ and $\omega$ since  $\nabla_xF$ is assumed to be boundedness-preserving.   The dominated convergence theorem for the stochastic integrals \cite[Ch.IV Theorem32]{protter} ensures that it converges in probability to $\int_{]0,t]} {\nabla}_xF_u(X_{u-},A_{u-}).dX(u)$. Using remark \ref{rightcontinuousito.remark} concludes the proof.\\
Consider now the general case.
Let $K_n$ be an increasing sequence of compact sets with $\bigcup_{n
\geq 0} K_n = U$ and denote
$$\tau_n=\inf\{s < t| X_s \notin K^n \mbox{ or } |A_s| >
n\}\wedge t$$ which are optional times. Applying the previous result
to the stopped process $(X_{t \wedge\tau_n},A_{t\wedge\tau_n})$
leads to: \ba F_t(X_{t \wedge\tau_n},A_{t\wedge\tau_n})-F_0(X_0,A_0)
= \int_{0}^{t\wedge\tau_n} \mathcal{D}_uF_u(X_u,A_u) du +
\frac{1}{2} \int_{0}^{t\wedge\tau_n}{\rm
tr}\left({}^t\nabla^2_xF_u(X_u,A_u) d[X](u)\right)   \nonumber \\
+ \int_{0}^{t\wedge\tau_n} {\nabla}_xF_u(X_u,A_u).dX+
\int_{t\wedge\tau^n}^t {\cal D}_uF(X_{u
\wedge\tau_n},A_{u\wedge\tau_n})du\ea The terms in the first line
converges almost surely to the integral up to time $t$ since
$t\wedge\tau_n=t$ almost surely for $n$ sufficiently large. For the
same reason the last term converges almost surely to 0.
\end{proof}


\appendix

\section{Some results on cadlag functions} \label{cadlagfunctions.sec}
For a cadlag function $f:[0,T]\mapsto \mathbb{R}^d$ we shall denote $\Delta f(t)=f(t)-f(t-)$ its discontinuity at $t$.
\begin{lemma}
For any cadlag function  $f:[0,T]\mapsto \mathbb{R}^d$
\begin{equation}\label{uniformcadlag.eq}
\nabla \epsilon > 0,\quad \exists \eta >0,\quad |x-y|\leq\eta \Rightarrow |f(x)-f(y)|\leq\epsilon+\sup_{t\in [x,y]}\{|\Delta f(t)|\}
\end{equation}\label{uniform.lemma}
\end{lemma}
\begin{proof} Assume the conclusion does not hold. Then there exists a sequence $(x_{n},y_{n})_{n\geq 1}$  such that $x_{n}\leq y_{n}$, $y_{n}-x_{n}\rightarrow 0$ but $|f(x_{n})-f(y_{n})|>\epsilon+\sup_{t\in [x_{n},y_{n}]}\{|\Delta f(t)|\}$.  We can  extract a convergent subsequence $(x_{\psi(n)})$ such that $x_{\psi(n)}\rightarrow x$. Noting that either an infinity of terms of the sequence are less than $x$ or an infinity are more than $x$, we can extract {\it monotone} subsequences  $(u_{n},v_{n})_{n\geq 1}$ of $(x_n,y_n)$ which converge to $x$.
 If $(u_n),(v_n)$ both converge to $x$ from above or from below, $|f(u_{n})-f(v_{n})|\rightarrow 0$ which yields a contradiction. If  one converges  from above and the other from below,
  $\sup_{t\in [u_{n},v_{n}]}\{|\Delta f(t)|\}>|\Delta f(x)|$ but $|f(u_{n})-f(v_{n})|\rightarrow |\Delta f(x)|$, which results in a contradiction as well. Therefore \eqref{uniformcadlag.eq} must hold.
\end{proof}

The following lemma is a  consequence of lemma \ref{uniform.lemma}:
\begin{lemma}[Uniform approximation of cadlag functions by step functions]\label{Approximation}\  \\
Let $h$ be a  cadlag function on $[0,T]$. If $(t^n_k)_{n\geq 0, k=0..n}$ is a sequence of subdivisions $0=t^{n}_{0}<t_{1}<...<t^{n}_{k_{n}}=t$  of $[0,T]$ such that:
$$\mathop{\sup}_{0 \leq i \leq k-1} |t^{n}_{i+1}-t^{n}_{i}| \rightarrow_{n\rightarrow \infty} 0\qquad\mathop{\sup}_{u \in [0,T]\setminus \{t^{n}_{0},...,t^{n}_{k_{n}}\}} |\Delta f(u) | \rightarrow_{n \rightarrow \infty} 0$$
then
\begin{equation}
\sup_{u \in [0,T]} |h(u)-\sum_{i=0}^{k_{n}-1} h(t_{i})1_{[t^{n}_{i},t^{n}_{i+1})}(u)+h(t^{n}_{k_{n}})1_{\{t^{n}_{k_{n}}\}}(u)| \rightarrow_{n\rightarrow \infty} 0
\end{equation}
\end{lemma}

\section{Proof of theorem \ref{Measurability.thm}} \label{measurabilityproof.sec}

\begin{lemma} \label{stoppingtime.lemma}
Consider the canonical space $\mathcal{U}_T$ endowed with the natural filtration of the canonical process $X(x,t)=x(t)$. Let $\alpha \in \mathbb{R}$ and $\sigma$ be an optional time. Then the following functional:
\begin{equation}
\tau(x)=\inf\{t > \sigma, \quad |x(t)-x(t-)|> \alpha\}
\end{equation}
is a stopping time.
\end{lemma}

\begin{proof}
We can write that:
\begin{equation}
\{\tau(x) \leq t\}=\bigcup_{q\in\mathbb{Q}\bigcap [0,t)}(\{\sigma \leq t-q\}\bigcap \{\sup_{t\in (t-q,t]}|x(u)-x(u-)|> \alpha\}
\end{equation}
and
\begin{equation}
\{\sup_{u\in (t-q,t]}|x(u)-x(u-)|>
\alpha\}=\bigcup_{n_0>1} \bigcap_{n>n_0} \{\sup_{1\leq
i\leq  2^n}|x(t-q\frac{i-1}{2^n})-x(t-q\frac{i}{2^n})|> \alpha\}
\end{equation}
thanks to the lemma \ref{uniform.lemma} in Appendix
\ref{cadlagfunctions.sec}.
\end{proof}
We can now prove Theorem \ref{Measurability.thm}  using lemma \ref{uniform.lemma} from Appendix \ref{cadlagfunctions.sec}.\\
{\bf Proof of Theorem \ref{Measurability.thm}:}
Let's first prove point 1.; by lemma \ref{cadlag.prop} it implies point 2. for right-continuous functionals and point 3. for left-continuous functionals. Introduce the following random subdivision of $[0,t]$:
\begin{eqnarray}
\tau^N_{0}(x,v)=0 \nonumber \\
\tau^N_{k}(x,v)=\inf\{t > \tau^N_{k-1}(x,v) | 2^{N}t \in \mathbb {N}\  {\rm or}\quad |v(t)-v(t-)|\vee |x(t)-x(t-)| > \frac{1}{N} \}\wedge t
\end{eqnarray}
From lemma \ref{stoppingtime.lemma}, those functionals are stopping times for the natural filtration of the canonical process.
We define the stepwise approximations of $x_{t}$ and $v_{t}$ along the subdivision of index $N$:
\begin{eqnarray}
x^{N}(s)=\sum_{k=0}^{\infty} x_{\tau^N_{k}(x,v)}1_{[\tau^N_{k}(x,v),\tau^N_{k+1}(x,v)[}(s)+x(t)1_{\{t\}}(s) \nonumber \\
v^{N}(s)=\sum_{k=0}^{\infty} v_{\tau^N_{k}(x,v)}1_{[\tau^N_{k}(x,v),\tau^N_{k+1}(x,v)[}(t)+v(t)1_{\{t\}}(s)
\end{eqnarray}
as well as their truncations of rank $K$:
\begin{eqnarray}
_{K}x^{N}(s)=\sum_{k=0}^{K} x_{\tau^N_{k}}1_{[\tau^N_{k},\tau^N_{k+1}[}(s) \nonumber \\
_{K}v^{N}(t)=\sum_{k=0}^{K} v_{\tau^N_{k}}1_{[\tau^N_{k},\tau^N_{k+1}[}(t)
\end{eqnarray}
First notice that:
\begin{equation}
F_t(x^{N}_{t},v^{N}_{t})=\lim_{K\rightarrow \infty} F_{t}(_Kx^{N}_{t},_Kv^{N}_{t})
\end{equation}
because $(_{K}x^{N}_{t},_{K}v^{N}_{t})$ coincides with
$(x^{N}_{t},v^{N}_{t})$ for $K$ sufficiently large. The truncations
$$F^n_{t}(_{K}x^{N}_{t},_{K}v^{N}_{t}) $$
are $\mathcal{F}_t$-measurable as they are
continuous functionals of the measurable functions:
$$\{(x(\tau^N_{k}(x,v)),v(\tau^N_{k}(x,v))),k \leq K\}$$
so their limit $F_{t}(x^{N}_{t},v^{N}_{t})$ is also $\mathcal{F}_t$-measurable. Thanks to lemma \ref{Approximation},
$x^{N}_{t}$ and $v^{N}_{t}$ converge uniformly to $x_{t}$ and
$v_t$, hence $F_{t}(x^{N}_{t},v^{N}_{t})$ converges to $F_t(x_{t},v_{t})$ since $F$ is continuous at fixed times.\\

Now to show optionality of $Y(t)$ for a left-continuous functional, we will exhibit it as limit of right-continuous adapted processes. For $t \in [0,T]$, define $i^n(t)$ to be the integer such that $t \in [\frac{iT}{n},\frac{(i+1)T}{n}).$ Define the process:
$Y^n((x,v),t)=F_{\frac{i^n(t)T}{n}}(x_{\frac{(i^n(t))T}{n}},v_{\frac{(i^n(t))T}{n}})$, which is piecewise-constant and has right-continuous trajectories, and is also adapted by the first part of the theorem. Now, by $d_{\infty}$ left-continuity of $F$, $Y^n(t) \rightarrow Y(t)$,
which proves that $Y$ is optional.\\
We similarly prove predictability of $Z(t)$ for a right-continuous functional. We will exhibit it as a limit of left-continuous adapted processes. For $t \in [0,T]$, define $i^n(t)$ to be the integer such that $t \in (\frac{iT}{n},\frac{(i+1)T}{n}]$. Define the process:
$Z^n((x,v),t)=F_{\frac{(i^n(t)+1)T}{n}}(x_{t-,\frac{(i^n(t)+1)T}{n}-t},v_{t-,\frac{(i^n(t)+1)T}{n}-t})$,
which has left-continuous trajectories since as $s \rightarrow t-$, $t-s$ sufficiently small, $i^n(s)=i^n(t)$ and $(x_{s-,\frac{(i^n(s)+1)T}{n}-s},v_{s-,\frac{(i^n(s)+1)T}{n}-s})$ converges to $(x_{t-,\frac{(i^n(t)+1)T}{n}-t},v_{t-,\frac{(i^n(t)+1)T}{n}-t})$ for $d_{\infty}$.
Moreover, $Z^n(t)$ is $\mathcal{F}_t$-measurable by the first part of the theorem, hence $Z^n(t)$ is predictable.
Since $F\in \mathbb{F}^\infty_r$, $Z^n(t) \rightarrow Z(t)$, which proves that $Y$ is predictable.\\

\section{Measure-theoretic lemmas used in the proof of theorem \ref{itonoprobability.theorem} and \ref{itojumpsnoprobability.theorem}} \label{measuretheory.sec}

\begin{lemma} \label{vagueconvergence.lemma}
Let $f$ be a bounded left-continuous function defined on $[0,T]$, and let $\mu(n)$ be a sequence of Radon measures on $[0,T]$ such that $\mu_n$ converges vaguely to a Radon measure $\mu$ with no atoms. Then for all $0 \leq s < t \leq T$,  with $\mathcal{I}$ being $[s,t]$, $(s,t]$ ,$[s,t)$ or $(s,t)$:
\ba \lim_n \int_\mathcal{I} f(u)d\mu_n(u) = \int_\mathcal{I} f(u)d\mu(u) \ea
\end{lemma}

\begin{proof}
Let $M$ be an upper bound for $|f|$, $F_n(t)=\mu_n([0,t])$ and
$F(t)=\mu([0,t])$ the  cumulative distribution functions associated
to $\mu_n$ and $\mu$. For $\epsilon > 0$ and  $u \in (s,t]$, define:
\ba \eta(u)=\inf \{{h > 0} ||f(u-h)-f(u)| \geq \epsilon \} \wedge u
\ea and we have $\eta(u) > 0$ by right-continuity of $f$. Define
similarly $\theta(u)$: \ba \theta(u)=\inf \{{h > 0} ||f(u-h)-f(u)|
\geq \frac{\epsilon}{2} \} \wedge u \ea By uniform continuity of $F$
on $[0,T]$ there also exists $\zeta(u)$ such that $\forall v \in
[T-\zeta(u),T], F(v+\zeta(u))-F(v) < \epsilon \eta(u)$. Take a
finite covering \ba [s,t] \subset \bigcup_{i=0}^{N}
(u_i-\theta(u_i),u_i+\zeta(u_i)) \ea where the $u_i$ are in $[s,t]$,
and in increasing order, and we can choose that $u_0=s$ and $u_N=t$.
Define the decreasing sequence $v_j$ as follow: $v_0 = t$, and when
$v_j$ has been constructed, choose the minimum index $i(j)$ such
that $v_j \in (u_{i(j)},u_{i(j)+1}]$, then either $u_{i(j)} \leq v_j
- \eta(v_j)$ and in this case $v_{j+1}=u_{i(j)}$, else $u_{i(j)} >
v_j-\eta(v_j)$, and in this case $v_{j+1}=\max(v_j-\eta(v_j),s)$.
Stop the procedure when you reach $s$, and denote $M$ the maximum
index of the $v_j$. Define the following piecewise constant
approximation of $f$ on $[s,t]$: \ba
g(u)=\sum_{j=0}^{M-1}f(v_j)1_{(v_{j+1},v_{j}]}(u) \ea Denote $J_1$
the set of indices $j$ where $v_{j+1}$ has been constructed as in
the first case, and $J_2$ its complementary. If $j \in J_1$,
$|f(u)-g(u)| < \epsilon$ on $[v_j-\eta(v_j),v_j]$, and
$v_j-\eta(u_{i(j)}))-v_{j+1} < \zeta(u_{i(j)+1})=\zeta(v_{j+1})$,
because of the remark that $v_j - \eta_{v_j} < u_{i(j)} -
\theta(u_{i(j)})$. Hence: \ba \int_{(v_j,v_{j+1}]}|f(u)-g(u)|d\mu(u)
\leq \epsilon [F(v_{j+1})-F(v_j)]+2M \epsilon \eta(v_{j+1}) \ea

If $j \in J_2$, $|f(u)-g(u)| < \epsilon$  on $[v_{j+1},v_{j}]$.  So
that summing up all terms we have the following inequality: \ba
\int_{[s,t]} |f(u)-g(u)|d\mu(u) \leq \epsilon \left(
F(t)-F(s)+2M(t-s)\right) \label{majorationvagueconvergence.eq} \ea
because of the fact that: $\eta(v_{j}) \leq v_{j}-v_{j+1}$ for $j <
M$. The same argument applied to $\mu_n$ yields:
\ba \int_{[s,t]} |f(u)-g(u)|d\mu_n(u) \leq \epsilon [F_n(t)-F_n(s-)] \nonumber \\
+2M \sum_{j=0}^{M-1} F_n(v_{j+1})-F_n(v_{j+1}-\zeta(v_{j+1}))  \ea
so that the $\limsup$ satisfies \eqref{majorationvagueconvergence.eq} since $F_n(u)$ converges to $F(u)$ for every $u$.

On other hand, it is immediately observed that \ba \lim_n
\int_{\mathcal{I}} g(u)d\mu_n(u) = \int_{\mathcal{I}} g(u)d\mu(u)
\ea since $F_n(u)$ and $F_n(u-)$ both converge to $F(u)$ since $\mu$
has no atoms ($g$ is a linear combination of indicators of
intervals). So the lemma is established.

\end{proof}

\begin{lemma} \label{vagueconvergence2.lemma}
Let $(f_n)_{n\geq 1}, f$ be left-continuous functions defined on $[0,T]$, satisfying:
\ba
\forall t \in [0,T], \lim_n f_n(t)=f(t) \qquad
\forall t \in [0,T], f_n(t) \leq K
\ea
Let also $\mu_n$ be a sequence of Radon measures on $[0,T]$ such that $\mu(n)$ converges vaguely to a Radon measure $\mu$ with no atoms. Then for all $0 \leq s < t \leq T$, with $\mathcal{I}$ being $[s,t]$, $(s,t]$ ,$[s,t)$ or $(s,t)$:
\ba \int_\mathcal{I} f_n(u)d\mu_n(u) \rightarrow_{n \rightarrow \infty} \int_s^t f(u)d\mu(u) \ea

\end{lemma}

\begin{proof}
Let $\epsilon > 0$ and let $n_0$ such that $\mu(\{\sup_{m \geq n_0}
|f_m-f| > \epsilon\}) < \epsilon$. The set $\{\sup_{m \geq n_0}
|f_m-f| > \epsilon\}$ is a countable union of disjoint intervals
since the functionals are left-continuous, hence it is a  continuity
set of $\mu$ since $\mu$ has no atoms; hence, since $\mu_n$
converges vaguely to $\mu$  \cite{billingsey99}: \ba \lim_n
\mu_n(\{\sup_{m \geq n_0} |f_m-f| > \epsilon\}) = \mu(\{\sup_{m \geq
n_0} |f_m-f| > \epsilon\}) < \epsilon \ea since $\mu_n$ converges
vaguely to $\mu$ which has no atoms.

So we have, for $n \geq n_0$: \ba \int_{\mathcal{I}}
|f_n(u)-f(u)|d\mu_n(u) \leq 2K \mu_n\left(\{\sup_{n \geq n_0}
|f_n-f|
> \epsilon\}\right) + \epsilon \mu_n(\mathcal{I}) \ea Hence the
$\limsup$ of this quantity is less or equal to: \ba 2K \mu(\{\sup_{m
\geq n_0} |f_m-f| > \epsilon\} + \epsilon \mu(\mathcal{I}) \leq (2K
+ \mu(\mathcal{I})) \epsilon  \ea On other hand: \ba \lim_n
\int_{\mathcal{I}} f(u) d\mu_n(u)=\int_{\mathcal{I}} f(u) d\mu(u)
\ea by application of lemma \ref{vagueconvergence.lemma}.

\end{proof}

\end{document}